\documentclass[twoside,12pt]{article} 
\usepackage{amsfonts} 
\usepackage{graphicx}
\usepackage{color}

\newtheorem{thm}{Theorem} 
 
\newtheorem{lemma}[thm]{Lemma} 
\newtheorem{prop}[thm]{Proposition}

\begin{document}

 \def\reals{{\mathbb R}}
 \def\ch{{\cal H}}
 \def\cA{{\cal A}}
 \def\cD{{\cal D}}
 \def\cK{{\cal K}}
 \def\cC{{\cal C}}
 \def\cN{{\cal N}}
 \def\cR{{\cal R}}
 \def\cS{{\cal S}}
 \def\cT{{\cal T}}
 \def\cV{{\cal V}}
 \def\ta{{\cal T}_{\subset}}
 \def\cI{{\cal I}}
 \def\bC{{\bf C}}
 \def\axis{{\bf A}}
 \def\fibr{{\bf H}}
 \def\ba{{\bf a}}
 \def\bb{{\bf b}}
 \def\bc{{\bf c}}
 \def\be{{\bf e}}
 \def\d{{\delta}} 
 \def\ci{{\circ}} 
 \def\e{{\epsilon}} 
 \def\l{{\lambda}} 
 \def\L{{\Lambda}} 
 \def\m{{\mu}} 
 \def\n{{\nu}} 
 \def\o{{\omega}} 
 \def\s{{\sigma}} 
 \def\v{{\varphi}} 
 \def\a{{\alpha}} 
 \def\b{{\beta}} 
 \def\p{{\partial}} 
 \def\r{{\rho}} 
 \def\ra{{\rightarrow}} 
 \def\lra{{\longrightarrow}} 
 \def\g{{\gamma}} 
 \def\D{{\Delta}} 
 \def\La{{\Leftarrow}} 
 \def\Ra{{\Rightarrow}} 
 \def\x{{\xi}} 
 \def\c{{\mathbb C}} 
 \def\z{{\mathbb Z}} 
 \def\2{{\mathbb Z_2}} 
 \def\q{{\mathbb Q}} 
 \def\t{{\tau}} 
 \def\u{{\upsilon}} 
 \def\th{{\theta}} 
 \def\la{{\leftarrow}} 
 \def\lla{{\longleftarrow}} 
 \def\da{{\downarrow}} 
 \def\ua{{\uparrow}} 
 \def\nwa{{\nwtarrow}} 
 \def\swa{{\swarrow}} 
 \def\nea{{\netarrow}} 
 \def\sea{{\searrow}} 
 \def\hla{{\hookleftarrow}} 
 \def\hra{{\hookrightarrow}} 
 \def\sl{{SL(2,\mathbb C)}} 
 \def\ps{{PSL(2,\mathbb C)}} 
 \def\qed{{\hfill$\diamondsuit$}} 
 \def\pf{{\noindent{\bf Proof.\hspace{2mm}}}} 
 \def\ni{{\noindent}} 
 \def\sm{{{\mbox{\tiny M}}}} 
 \def\sc{{{\mbox{\tiny C}}}} 
 \def\ke{{\mbox{ker}(H_1(\p M;\2)\ra H_1(M;\2))}} 
 \def\et{{\mbox{\hspace{1.5mm}}}}

 \definecolor{Turquoise}{cmyk}{.5,0,0,.5}
 \definecolor{Red}{cmyk}{0,.5,.5,.5}
 \definecolor{LRed}{cmyk}{0,1,1,0}
 \definecolor{LGreen}{cmyk}{1,0,1,0}
 \definecolor{LBlue}{cmyk}{1.0,1.0,0,0}
 \definecolor{LMagenta}{cmyk}{0.2,0.4,0,0.2} 
 \definecolor{DMagenta}{cmyk}{0.16,0.44,0,.42}
 \definecolor{LOrange}{cmyk}{0,.2,.4,0}

\pagestyle{myheadings}
\markboth{William W. Menasco \today}{An addendum on iterated torus knots \today}
\title{An addendum on iterated torus knots} 
 
\vspace{1mm} 
 
\author{William W. Menasco \thanks{partially supported by NSF grant \#DMS 0306062} \\  
{\small University at Buffalo} 
 \\ {\small Buffalo, New York 14260} }
\date{\today}

\maketitle
 
\section{Introduction}
\label{section: Introduction}
In Theorem 1.2 of the paper \cite{[M1]} the author claimed to have proved that all transversal knots whose topological knot type is that of an iterated torus knot (we call them {\it cable knots}) are transversally simple.  That theorem is false, and the Erratum  \cite{[M2]} identifies the gap.  The purpose of this paper is to explore the situation more deeply, in order to pinpoint exactly which cable knots are {\it not} transversally simple.   The class is subtle and interesting.   We will recover the strength of the main theorem in \cite{[M1]}, in the sense that we will be able to prove a strong theorem about cable knots, but the theorem itself is more subtle than Theorem 1.2 of \cite{[M1]}.

\begin{figure}[!ht]  
\centerline{\includegraphics[scale=.9, bb= 0 0 451 210]{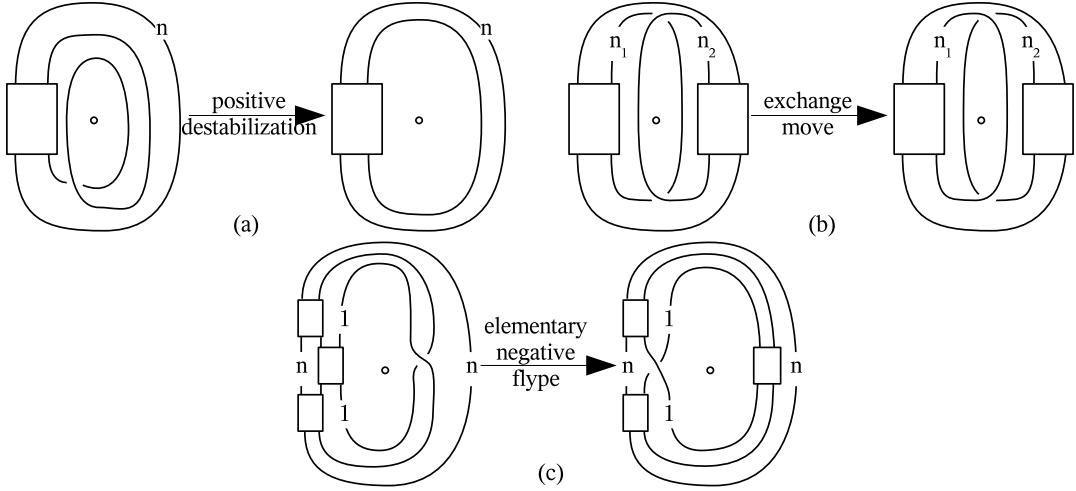}}
\caption{\footnotesize The illustration indicates that some of the strands in the
block-strand diagrams can be weight, i.e. multiple parallel strands.}
\label{figure:isotopies}
\end{figure}

Before we can state the main results in this paper we need to set up some machinery. 
The reader should recall that a closed braid $X$ is {\em exchange reducible} to a closed braid
$Y$ if through a sequence of exchange moves, destabilizations and braid isotopies we can
take $X$ to $Y$.  In the case of cable knots saying that a specified class
is {\em exchange reducible} means that it is exchange reducible to
the standard braid representation of minimal index.  See Figure \ref{figure:isotopies}.
Recall also that in \cite{[BW]} it was established that exchange reducible
closed braids representing oriented knots are transversally simple when the braids are
viewed as transversal knots in the standard contact structure of $\reals^3$ or $S^3$.
The key question for us is to understand which cable knots are exchange-reducible and which ones are not.   To specify the large subclass of cable knots that are not exchange reducible we need to
develop a new mode of representing knots.
It is first useful to set up a coordinate system for
$ (\reals^3, \{z{\rm -axis}\}) \subset (S^3, \axis)$.  Let $(\rho , \theta , z)$ be the
classical cylindrical coordinate system for $\reals^3$.  We specify notation
for the unit cylinder, $\cT_1 = \{(\rho , \theta , z) | \rho = 1 \}$; and level
planes, $P_{z_0} = \{ (\rho , \theta , z) | z = z_0(=\{{\rm constant}\}) \}$.

As in \cite{[M1]} $\axis ( = z{\rm -axis} \cup \{\infty\})$ will be notation for our
braid axis, and $\{ H_\theta | 0 \leq \theta < 2 \pi \}$ will be our collection of disc
that make up $\fibr$, our braid fibration.

For a closed $n$-braid $X \subset \reals^3 \setminus \{ z{\rm -axis} \}$ a
{\em rectangular block presentation} $\cR_X$ is a collection of disc
$\{R_1 , \cdots , R_l \} \subset \reals^3 \setminus \{ z{\rm -axis} \}$ satisfying the
following conditions.

\begin{itemize}
\item[(0)] For every point $(\rho , \theta , z) \subset R_i (\in \cR_X)$ we have $\rho \geq 1$. 
\item[(1)] For each $R_i$ we have $\partial R_i = r_i^1 \cup r_i^2 \cup r_i^3 \cup r_i^4$
where:
\begin{itemize}
\item[(a)] $r_i^1 \subset H_{\theta_i^1}$ and $r_i^3 \subset H_{\theta_i^3}$ for
$H_{\theta_i^1} , H_{\theta_i^3} \in \fibr$ with $\theta_i^1 < \theta_i^3$.  (That is,
$r_i^1$ is the {\em bottom side} and $r_i^3$ is the {\em top side} of $R_i$ in the $\fibr$
fibration.)
\item[(b)] $r_i^2 \subset P_{z_i^2} \cap \cT_1$ and $r_i^4 \subset P_{z_i^4} \cap \cT_1$.
\end{itemize}
\item[(2)] Still appealing to the specified notation in condition (1),
$\theta_i^1 = \theta_j^1$ iff $i=j$.  Likewise, $\theta_i^3 = \theta_j^3$ iff
$i = j$.  Similarly, $z_i^2 = z_j^2$ iff $i=j$; and, $z_i^4 = z_j^4$ iff $i=j$. 
\item[(3)] There exists a fixed integer $0 \leq k \leq l$ such that if
$\theta_i^4 = \theta_j^1$ then $i+k=j \ {\rm mod}\ l$.  In particular,
$r_i^4 , r_{i+k}^1 \subset H_{\theta_i^4} (=H_{\theta_{i+k}^1})$.
\item[(4)] $r_i^2 \cap r_{j}^4 \not= \emptyset$ then $i = j+1$.  In particular,
$z_{i+1}^2 = z_i^4 \ {\rm mod} \ l$.
\item[(5)] The leaves of the induced foliation of each $R_i$ are arcs parallel to
the $r_i^1$ side (or, equivalently the $r_i^3$ side).  Thus, each leaf
of $\cR_X = \cup_{1 \leq i \leq l} R_i$ is also an arc.
\item[(6)] By a braid isotopy of $X$ in $\reals^3 \setminus \{z {\rm -axis} \}$ we can position
it so as to have $X \subset \cup_{1 \leq i \leq l} R_i ( = \cR_X)$ such that $X$ transversely
intersects each leaf in the foliation of $\cR_X$ exactly once.  (Thus,
$X$ is a deformation retract of $\cR_X$, the retraction occurring along leaves of
$\cR_X$ foliation.)
\end{itemize}

Rectangular block presentations of any given closed $n$-braid are readily produced.
Referring ahead Figure \ref{figure:discs plus blocks} if the reader only focuses on the
rectangular blue blocks then one has a local portion of a possible rectangular
block presentation with $k=1$.  More specifically, again just focusing on the blue
blocks, in Figure \ref{figure:trefoil torus} one has a rectangular block presentation of the
positive trefoil.

We say that a rectangular block presentation $\cR_X$ has {\em homogeneous twisting}
if for any triple $(z_{i-1}^4, z_i^4, z_{i+1}^4)$ satisfying the condition that
$z_{i+1}^4 < z_{i-1}^4 < z_i^4$ then we have that $ z_{i+1}^4 < z_{i+2}^4 < z_i^4$.
(Similarly, if $z_i^4 < z_{i-1}^4 < z_{i+1}^4$ then we have that $ z_i^4 < z_{i+2}^4 < z_{i+1}^4$.)
In \S\ref{section: building intlocking presentations} we show how rectangular block presentations having homogeneous twisting of any
$n$-braid are readily produced.  The rectangular block presentation of the positive trefoil
in Figure \ref{figure:trefoil torus} is such a presentation.

Now we introduce a concept that is more restrictive on rectangular block presentations.
Given $\cR_X$, as we traverse the copy of $X$ it contains we will cross the
leaves that contain the $r_i^1$ and $r_i^3$ sides of $R_i {\rm 's}$.  Our notation for this
finite sequence of arc leaves will be $\{ \l_1 , \cdots , \l_l \} \subset \cR_X$ (where
the index indicates the cyclic order of intersection with $X$).  By our
conditions on $\cR_X$, $\l_i$ will: intersect $k+1$ blocks; contain $r_i^4$ and $r_{i+k}^1$
(or, flipping things around, $r_i^1$ and $r_{i+k}^4$; and, $\l_i \cap \cT_1$ will
be $3+k$ points $\{ \varrho_i(1) , \cdots , \varrho_i(3+k) \}$.
    
We say that a rectangular block presentation $\cR_X$ is {\em interlocking} if for
every consecutive pair of leaves $(\l_j , \l_{j+1})$ there exists a sub-collection
of leaves $\{ \l_{u_j^1} , \cdots , \l_{u_j^v} \}$ such that:
\begin{itemize}
\item[(i)] For disc fibers of $\fibr$, $\l_j \subset H_{\theta_j}, \ \l_{u_j^1} \subset H_{\theta_j^1}, \cdots , \ \l_{u_j^v} \subset H_{\theta_j^v}, \ \l_{j+1} \subset H_{\theta_{j+1}}$,
we have that the angular order being
$\theta_j < \theta_j^1 < \cdots < \theta_j^v < \theta_{j+1}$.
\item[(ii)] The endpoint $\varrho_{u_j^1}(3+k) \in \l_{u_j^1}$ has $z$-coordinate between
the endpoints $\varrho_j(1)$ \& $\varrho_j(2)$ of $\l_j$.
\item[(iii)] Iteratively,
the endpoint $\varrho_{u_j^{w+1}}(3+k) \in \l_{u_j^{w+1}}$ has $z$-coordinate between
the endpoints $\varrho_{u_j^{w}}(1)$ \& $\varrho_{u_j^{w}}(2)$ of $\l_{u_j^{w}}$.
\item[(iv)] Finally,
the endpoint $\varrho_{u_j^v}(3+k) \in \l_{u_j^v}$ has $z$-coordinate between
the endpoints $\varrho_{j+1}(3+k -1)$ \& $\varrho_{j+1}(3+k)$ of $\l_{j+1}$.
\end{itemize}

Not all closed $n$-braids have interlocking rectangular block presentations.
For example, we state without proof the the unknot does not have such a presentation.
(Once the naturality of interlocking is understood this claim can be see as a corollary of
Theorem 1 \cite{[M2]}.)  Again, focusing only on the
blue rectangular blocks in Figure \ref{figure:trefoil torus}, it is readily
observed that the rectangular block presentation of the positive trefoil has homogeneous
twisting and is interlocking.  (The intuition behind the concept of an interlocking presentation should
be self-evident: if we try to push $r_j^3$ forward in $\fibr$ to slide past $r_{j+1}^1$ we are
forced to push $\l_{u_j^1}$ forward, which pushes $\l_{u_j^2}$, which pushes, etc; until
we are forced to push $\l_{j+1}$ forward.)

We can now state the main result in this paper.

\begin{thm}
\label{theorem:addendum theorem}
Every cable knot is either exchange reducible or through
a sequence of braid isotopies, exchange moves and $\pm$-destabilizations reducible
to an iterated cabling of a braid
that admits an interlocking homogeneous twisting rectangular block presentation.
(This secondary braid is necessarily also a cable knot.)  Torus knots, in particular, are exchange reducible.
\end{thm}

The ending observation, that torus knots are exchange reducible, immediately implies the
result of \cite{[M2]},
torus knots are transversally simple. 

To state this result in a more exacting manner, it is useful to recall our notational machinery
from \cite{[M1]}. Let $\cT_{\cC} \subset S^3$ be a
peripheral torus for an oriented knot $\cC \subset S^3$.
The oriented simple closed curve on $\cT_{\cC}$
that represents the homotopy class of $pm+ql$, where $m$ is the meridian homotopy class,
$l$ is the preferred longitude homotopy class and $p,q \in \z$,
is called the {\em $(p,q)$ cable of $\cC$}.
(Since the closed curve is necessarily of one component, $p$ and $q$ are relatively prime.)
We used the notion $\bC(\cC,(p,q))$ to
indicate the resulting oriented knot of this {\em cabling operation}.

Now, an iterated torus knot is obtained
by taking the initial knot $\cC_0$ as the oriented unknot.
Specifically, we take sequence of co-prime $2$-tuples of integers
$(P,Q) = \{(p_1,q_1),(p_2,q_2),\cdots,(p_h,q_h)\}$, and
we can construct the oriented knot
$$ K_{(P,Q)} = \bC(\bC(\cdots\bC(\bC(\cC_0,(p_1,q_1)),(p_2,q_2))\cdots,(p_{h-1},q_{h-1})),
(p_h,q_h)).$$

Theorem \ref{theorem:addendum theorem} can then be view as saying any $K_{(P,Q)}$ is
either exchange reducible or there is a proper subsequence of $(P,Q)$,
$(P^\prime,Q^\prime) = \{(p_1,q_1), \cdots , (p_{h^\prime},q_{h^\prime})\}$ such that
$K_{(P^\prime,Q^\prime)}$ admits an interlocking homogeneous twisting rectangular block presentation.

The outline of this note is as follows.  In \S\ref{section:mistake} we detail the
error in the analysis of standardly tiled cabling tori that lead
to the flaw in the proof of Theorem 1.1 of \cite{[M2]}.  Once there is a sufficient understanding of the
corrupting error we go on to prove that all torus knots are exchange reducible.
In \S\ref{section:addendum analysis} we continue our
analysis of the standardly tiled
tori that carry $K_{(P,Q)}$.  This will allow us to then give the proof
of Theorem \ref{theorem:addendum theorem}.

In \S\ref{section: building intlocking presentations}
we give a procedure for constructing interlocking homogeneous twisting rectangular
block presentations.  Such presentations have applications in the study of transversal and
Legendrian knots in the standard contact structure for $\reals^3$ or $S^3$.
In particular, in the \S\ref{section: The Etnyre-Honda cable.}-Appendix (joint with H. Matsuda)
we will discuss the explicit representations to the $(2,3)$ cabling of the
$(2,3)$ torus knot (the positive trefoil) that was implicitly discovered in
Theorem 1.7 of \cite{[EH]}.
The utilized rectangular block presentation was first discovered V. Pinciu \cite{[P]} and are
further exploited in \cite{[C],[D]}.

We remark that with the availability of manuscripts on-line it makes sense to make use
of color in figures that can be viewed on a computer screen in spite of the fact that
color is not as widely available in print.  Thus, color labeling the different salient features
of a number of figures will be utilized in this note.  Where possible we will also employ
the redundancy of symbolic labeling. 

%

\section{Understanding the error in Theorem 1.1 and exchange reducible for torus knots}
\label{section:mistake}
Before we describe the error in Theorem 1.1 of \cite{[M2]} we briefly review our notation.
The triple $(K,\ta,m)$ is comprised of an oriented cable knot $K$ that is a curve
on the torus $\ta \subset S^3$ having a meridian meridian curve $m$.  The quadruple
$(K,\ta,m,\Delta_m)$ is a triple with the addition of $2$-disc $\Delta_m$ satisfying
$ m = \ta \cap \Delta_m$.  When $\bb \bb $-tiles are the only
tiles in the singular foliation on $\ta$, the tiling determines the four graphs
$G_{\delta , \epsilon}$, $ (\delta ,\epsilon) \in \{ (+,+), (+,-), (-,+), (-,-) \} $.
Moreover, when $\ta$ has a $\bb\bb$-tiling the $\cS_K$ denotes the $\bb$-support of $K$ in $\ta$,
the closure of the union of all the $\bb$-arcs that $K$ intersects.  Similarly, $S_m$ denotes
the $\bb$-support of the meridian curve $m$.

The error occurs in Lemma 5.13\cite{[M1]}.  In particular, the alteration from Figure 15b\cite{[M1]} to 16b\cite{[M1]}
is not valid due to an obstruction which we will describe shortly.  Picking up the argument
at Lemma 5.13\cite{[M1]}, we have established in Propositions 5.5\cite{[M1]} through 5.10\cite{[M1]} and Lemma 5.11\cite{[M1]} that:
\begin{itemize}
\item All of the components of our graphs $G_{\delta , \epsilon }$ are homeomorphic to
$S^1$ on $\ta$.
\item $K$ coherently intersects components of $G_{\d,\e}$.
\item $\cS_K$ is either an annulus having one boundary component from $G_{\d,\e}$ and one
from $G_{-\d,-\e}$; or it is a torus-minus-a-disc having its single boundary curve being a
union of four arcs, one arc coming from each of the four graphs $G_{\d,\e}$.
\item The meridian curve $m$ intersects each component of the four graphs $G_{\d,\e}$
and intersects them coherently.
\item The intersection $\cS_K \cap \cS_m$ is a disjoint union of $\bb$-rectangles, rectangular
regions in the foliation that is the closure of the union of all the $\bb \bb$-arcs in a single
homotopic family.
\end{itemize}
\begin{figure}[htb]
\centerline{\includegraphics[scale=.55, bb=0 0 831 596]{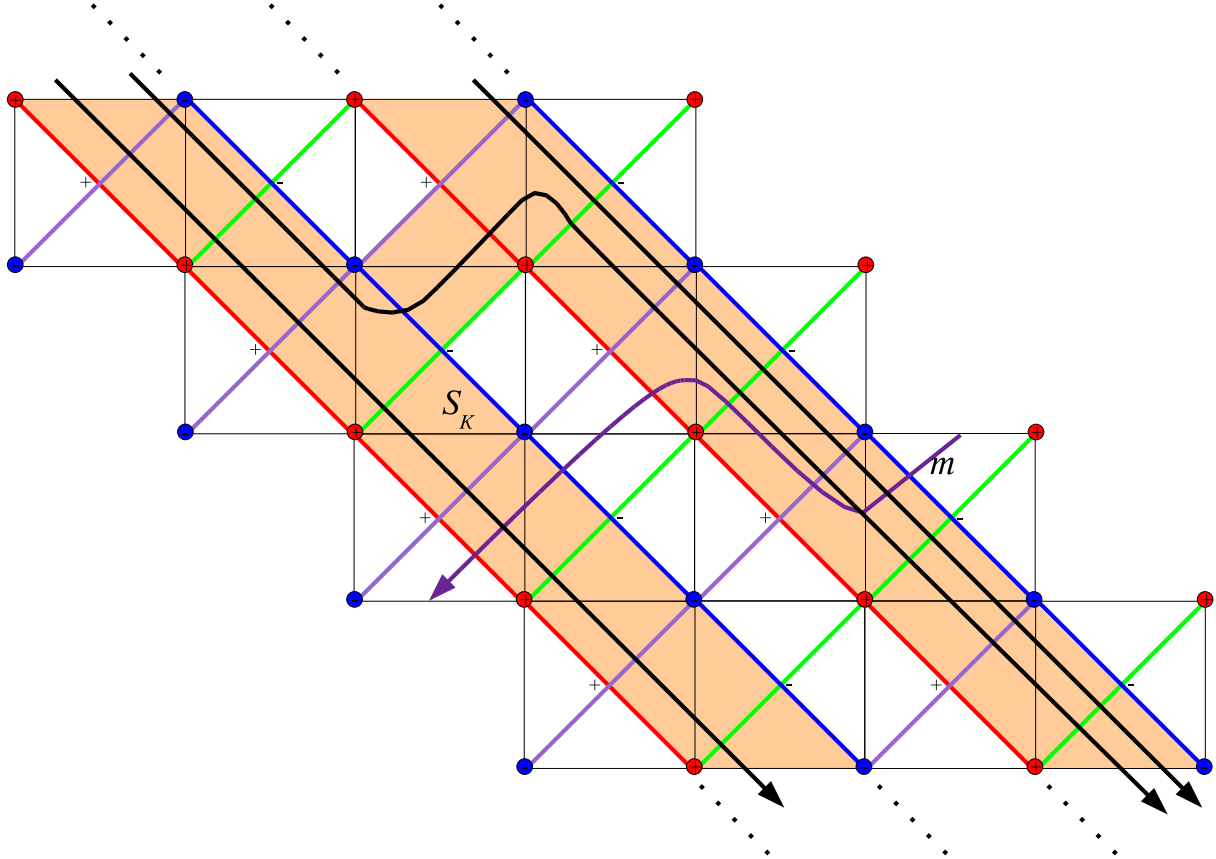}}
\caption{\footnotesize We depict the staircase pattern for the foliation of $\ta$ and
color code the four graphs of its foliation:  positive vertices are labeled $+$ and are
coded \textcolor{LRed}{red}; negative vertices are labeled $-$ and coded \textcolor{LBlue}{blue};
edges in $G_{+,+}$ are coded \textcolor{LRed}{red}; edges in $G_{-,-}$ are coded
\textcolor{LBlue}{blue}; edges in $G_{-,+}$ are coded \textcolor{LMagenta}{light magenta};
and edges in $G_{+,-}$ are coded \textcolor{LGreen}{green}.  We will always depict
the meridian curves as \textcolor{DMagenta}{dark magenta} and $X$ as
\textcolor{black}{black}.  Thus, are color code for the graphs is:
\textcolor{LRed}{$G_{+,+}$}; \textcolor{LBlue}{$G_{-,-}$}; \textcolor{LGreen}{$G_{+,-}$};
and \textcolor{LMagenta}{$G_{-,+}$}.  The $\bb$-support, $\cS_K$ is colored
\textcolor{LOrange}{light orange}.}
\label{figure:staircase}
\end{figure}

To help with the visualization of the foliation of $\ta$ and $\cS_K$ we refer to the result
in \cite{[N]} from which we know that the foliation $\ta$ has a {\em staircase pattern}
as illustrated in Figure \ref{figure:staircase}.  Specifically, by cutting
open the foliation of $\ta$ along two edge-path the are the union of $\bb$-arcs, one
resembling the steps of a staircase and the other a straight path, we can place flat on a plane
a fundamental region of $\ta$.

\begin{figure}[htpb]
\centerline{\includegraphics[scale=.7, bb=0 0 547 464]{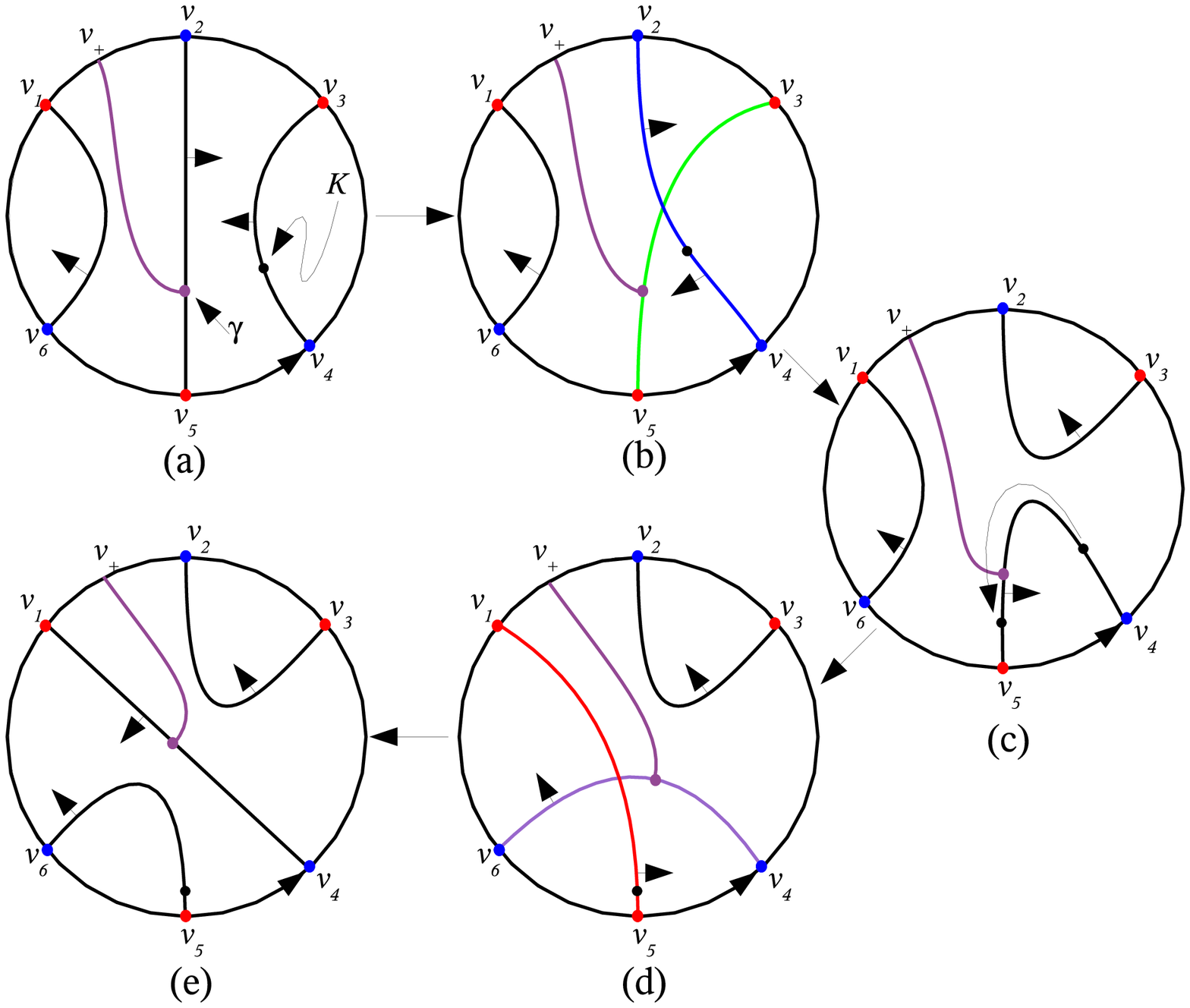}}
\caption{\footnotesize In (b) and (d) we have negative and positive singular leaves,
respectively.  Between the singularities in (b) and (d) the $\ba$-arcs adjacent to $v_+$
``block'' any other singularity occurrence from obstructing a change of fibration that
would allow us to first do the positive singularity first then the negative singularity.}
\label{figure:first Ht-sequence-non-standard}
\end{figure}

The proof of Lemma 5.13\cite{[M1]} employed the use of ``non-standard''
change of fibration that is illustrated in Figures 15 \& 16\cite{[M1]}.  Specifically, an arc
$\gamma \subset \partial S_m$ that in the foliation of $S_m$ is adjacent to a single
vertex (depicted as $v_+$ in Figure 15a\cite{[M1]}) and crosses in succession a negative singular
leaf followed by a positive singular leaf (depicted as $s_{i+1}$ in Figure 15b\cite{[M1]}).
It is useful to complement the illustration in Figure 15\cite{[M1]} with the corresponding
$H_\theta$-sequence.  In Figure \ref{figure:first Ht-sequence-non-standard} we illustrate
this sequence.  The salient feature of this sequence is that the $\ba$-arcs adjacent to
$v_+$ and having endpoints on $\gamma$ ``block'' the occurrence of any singularity between
the disc fibers in Figure \ref{figure:first Ht-sequence-non-standard}b and
\ref{figure:first Ht-sequence-non-standard}d that would act as an obstruction to
reversing to order of occurrence of the singularities in
\ref{figure:first Ht-sequence-non-standard}b and
\ref{figure:first Ht-sequence-non-standard}d.  Thus, this portion of
our erroneous argument was correct.  (When we more throughly analyze the geometry of
$\ta$ in \S\ref{section:addendum analysis} we will refer back to
Figure \ref{figure:first Ht-sequence-non-standard}.)

\begin{figure}[htpb]
\centerline{\includegraphics[scale=.6, bb=0 0 554 190]{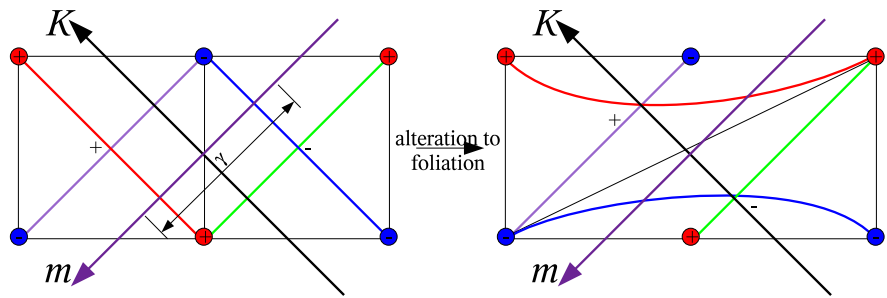}}
\caption{\footnotesize The left illustration depicts the knot $K$ and the portion
of the meridian curve $\g \subset m$ intersecting essentially.  Thus, when
we alter the foliation as depicted in Figure 16b\cite{[M1]} $K$ is no longer transverse to
the leaves of the foliation of $\ta$.}
\label{figure: obstruction}
\end{figure}

However, there is an obstruction that can possible prevent the alteration from
Figure 15b to 16b\cite{[M1]}
comes from the placement of
$K$ in the foliation of $\ta$.  As illustrated in Figure \ref{figure: obstruction},
when $K$ is positioned in the foliation of $\ta$ so that it transversely intersects
$m$ in an ``essential'' manner between the singularities $s_i$ and $s_{i+1}$,
after the alteration to Figure 16b\cite{[M1]} is performed
$K$ may no longer be transverse to the foliation and, thus, no longer positioned
as a braid in the fibration $\fibr$.  In other words, the arc $\gamma \subset m$ which we 
would hope to use for a non-standard change of fibration is invalid since it also
intersects $\bb$-support of $K$ in $\ta$.

The question then becomes when is such an obstruction unavoidable?  To answer this
we revisit Proposition 5.8\cite{[M1]} which governs the behavior of $\cS_K$, the $\bb$-support
of $K$.  By Proposition 5.8\cite{[M1]} $\cS_K$ is topologically either an annulus or a torus
minus an open disc.

When $\cS_K$ is topologically an annulus then, by Proposition 5.8(i)\cite{[M1]},
$K$ is in fact parallel in $\ta$ to a component of $G_{\e, \d}$ for some $\e,\d$ pair.
(Equivalently, $K$ will also be parallel to $G_{-\e,-\d}$.)
Thus. for an obstruction to the alteration in the foliation as depicted in
Figure 16b\cite{[M1]} to be unavoidable the arc $\g$ in
Figure 15b\cite{[M1]} must always be parallel in
the foliation of $\ta$ to a component of $G_{-\e ,\d}$ (or, equivalently, $G_{\e, -\d}$).

To see when $\g$ is always parallel to a component of $G_{-\e ,\d}$ we start with
the valence one vertex depicted in Figure 16a\cite{[M1]} and work backwards.
In Figure 16a\cite{[M1]} comes
from the existence of a valence two vertex near the meridian boundary curve $m$
that is adjacent to two singularities of common parity as depicted in Figure 15a\cite{[M1]}.  But,
by the discussion associated with Figure 14\cite{[M1]} such valence two vertices are forced to
exist.  Briefly summarizing that discussion, for $\Delta_m$ to be geometrically
realizable $m$ must transversally intersect $G_{\e, \d}$ for all four possible
pairs of $\e,\d$.  Moreover, in cyclic sequence $m$ intersects in an alternating
fashion $G_{+,+}$ and $G_{-,-}$; then $m$ intersects in an alternating
fashion $G_{+,-}$ and $G_{-,+}$.  Now if $m$ intersects any $G_{+,+}$ and $G_{+,-}$
more than once (equivalently, $G_{-,-}$ and $G_{-,+}$) then there will exist
valence two vertices adjacent to singularities of either parity; and, thus, we
will be able to choose $\g$ arcs parallel to the desired $G_{\e,\d}$ graph so as to
avoid an obstruction to the alteration in the foliation.  We conclude that the
obstruction is unavoidable exactly when $m$ intersects either $G_{+,+}$ and
$G_{-,-}$ once; or $G_{+,-}$ and $G_{-,+}$ once.
(We will refer back to this as {\em observation-$\star$}.)
Moreover, since every possible
$\g \subset m$ must encounter $\cS_K$ so as to create an obstruction we can
conclude that $\partial \cS_K = G_{\e,\d} \cup G_{-\e,-\d}$ for some $\e,\d$ pair,
i.e. $G_{\e,\d}$ and $G_{-\e,-\d}$ have exactly one component which is
topologically an $S^1$.

In the case where $\cS_K$ is a torus minus an open disc we have a similar analysis.
To start it is useful to notice that Proposition 5.8(ii)\cite{[M1]} implies the $\cS_K$ contains
a component of $\ta \setminus (G_{\e , \d} \cup G_{-\e ,-\d})$ for some $\e,\d$.
Next, to insure that all possible
$\g \subset m$ encounter an obstruction (as in the case of $\cS_K$ being an annulus)
for this designated $\e,\d$ pair, $G_{\e,\d}$ and $G_{-\e,-\d}$ both must be
topologically $S^1$ and, thus,
$\ta \setminus (G_{\e , \d} \cup G_{-\e ,-\d})$ can only have two components.
We conclude again that the
obstruction is unavoidable exactly when $m$ intersects either $G_{+,+}$ and
$G_{-,-}$ once; or $G_{+,-}$ and $G_{-,+}$ once.  Moreover, $\cS_K$ intersect
all outer-most $\bb$-arcs of $\ta$ by Lemma 2.2\cite{[M1]}.
(Again, we will refer back to this as {\em observation-$\star$}.)

We will continue this analysis in \S\ref{section:addendum analysis}.  But, we
are now in a position to establish that torus knots are exchange reducible, the
concluding statement of Theorem \ref{theorem:addendum theorem}.

%

\subsection{Torus knots are exchange reducible}
\label{subsection:torus knots}
We begin the argument with a brief review of the flow of the argument in \cite{[M1]}.
The cases of when the torus $\ta$ has a circular or mixed foliation are still
established, since error in Lemma 5.13\cite{[M1]} did not disturb these results.
The remaining case is still when $\ta$ has tiled foliation.

As analyzed in \S\ref{section:mistake}, there can exist an avoidable obstruction to performing
the alternation to the Figure 16\cite{[M1]}, since $K$ is by assumption a torus knot
$\ta \subset S^3$ is an unknotted torus.  Thus, we have available ``meridian disc'' on
both sides of $\ta$ in $S^3$.  Let $m^i \subset \ta$ correspond
to the {\em inner meridian curve} of \S5.2\cite{[M1]}; and designate $m^o \subset \ta$ as the
{\em outer meridian curve}.  (It may be helpful to think of $m^i$ as a meridian
curve of the unknotted core of the solid torus $\ta$ bounds, and $m^o$ is a
preferred longitude of the unknotted core.)  When viewing $m^o$ in the foliation
of $\ta$ we can apply the conclusions of Proposition 5.5, 5.8, 5.9 \& 5.10\cite{[M1]}.

Continuing let $\Delta_{m^o}$ be a disc that $m^o$ bounds.  Our tactic
is to now use $\Delta_{m^o}$ to produce a non-obstructed Figure 16\cite{[M1]} type
alteration to the foliation of $\ta$.  With this in mind we adapt Lemma 5.11\cite{[M1]}
to control the foliation of $\Delta_{m^o}$.  Thus, we have the following
lemma.

\begin{lemma}
\label{lemma:the intial foliation of Delta_m^o}
Let $(K,\ta,m,\Delta_{m^o})$ be a quadruple where the triple $(K,\ta,m^o)$ 
satisfies the conclusions of Propositions 5.5,
5.8, 5.9 and 5.10\cite{[M1]}.  Assume there is a $\e,\d$ pairing such
that the inner meridian curve $m^i \subset \ta$
intersects components of $G_{\e,\d}$ \& $G_{-\e,-\d}$ some number of times but only intersects $G_{-\e,\d}$ \& $G_{\e,-\d}$ once each.
We can then assume that the initial foliation of the spanning disc $\Delta_m^o$
satisfies the following conditions:
\begin{itemize}
\item[(a)] There are only $\ba\bb$- \& $\bb\bb$-singularities.
\item[(b)] The loop $l \subset \Delta_{m^o}$ that is parallel to $m^o$ and contains
all of the positive vertices adjacent to $m^o$ plus all of the $\ba\bb$-singularities
is the union of an edge-path in
$G_{+,+} \subset \Delta_m^o$ and an edge-path in $G_{+,-} \subset \Delta_m$.
\item[(c)] Either $G_{-,-} \cap \Delta_{m^o}$ or $G_{-,+} \cap \Delta_{m^o}$ contains
an edge-path that starts and end on $m^o$.
\item[(d)] The outer meridian curve $m^0 \subset \ta$ intersects components of
$G_{-\e,\d}$ \& $G_{\e,-\d}$ some number of times but only intersects $G_{\e,\d}$
\& $G_{-\e,-\d}$ once.
\end{itemize}
\end{lemma}

\pf
The argument is essentially a repeat of the argument in Lemma 4.11\cite{[M1]}.  Given
any $\Delta_{m^o}$ we can construct a new spanning disc $\Delta^\prime_{m^o}$
by extending $\Delta_{m^o}$ along an annulus as depicted in Figure 13\cite{[M1]}.  (To
briefly reiterate the argument concerning the existence of such a disc,
if $m^o$ were parallel to an component of $G_{\e,\d}$ then we will be able to
construct a tiled embedded disc in $S^3$ that has a loop in one of its
$G_{\e,\d}$ graphs, a contradiction of Lemma 3.8 of \cite{[BF]}.)
This repeated argument gives us statements (a) \& (b).

To get statement (c), assume without loss of generality, the loop $l$ of statement (b)
has at least two negative singularities.  Then let $C_{-,-} \subset \Delta_{m^o}$ be a component of the graph $G_{-,-}$.
Notice that $C_{-,-}$ must have an endpoint on $m^o$ (otherwise we again get a contradiction to
Lemma 3.8 of \cite{[BF]}).  If $C_{-,-} \cap m^o$ contains more than one point then $C_{-,-}$ will contain
an edge-path that satisfies statement (c).
So suppose that $C_{-,-} \cap \Delta^\prime_{m^o}$ is exactly
one point.  Let $\cN(C_{-,-})$ be a regular neighborhood of $C_{-,-}$ that is the closure of
the union of all of the $\ba$ and $\bb$ arcs that are adjacent to the vertices contained
in $C_{-,-}$.  Then $\partial \cN(C_{-,-})$ contains a component $C_{+,+} \subset G_{+,+}$ that contains both positive vertices which are adjacent to the unique $\ba\bb$-singular
leaf that by assumption $l$ intersects.  Now if we take $\Delta_{m^o}$ and repeat
the construction that extends it by the addition of a Figure 13\cite{[M1]} annulus, this $C_{+,+}$
component will be extended to a component $C_{-,+} \subset G_{-,+}$ that has
$C_{-,+} \cap m^o$ containing at least two points.  (By abuse of notation we
still refer to the resulting boundary curve as $m^o$.)  Thus, statement (c) is
established.

Finally, we establish statement (d) as follows.  Suppose a portion of
$m^o (= \partial \Delta_{m^o})$ intersects in sequence components of
$G_{\e,\d}$ and \& $G_{-\e,-\d}$ more than once.  Then by the
previous construction we can produce an edge-path $C$ in the foliation of
$\Delta_{m^o}$ that is either in the graph
$G_{-,-}$ or $G_{-,+}$; and begins and ends on $m^o$.  (For convenience of expository
let us assume that this edge-path is in $G_{-,-}$.)  But, by construction the
endpoints of such an edge-path can be positioned so as to be contained in
an arc $\g \subset m^i (\subset \Delta_{m^i})$ as depicted in Figure 15a and 16a\cite{[M1]}.
The two endpoints of this edge-path in $\Delta_{m^o}$ could then be ``coned'' to
the valence-one vertex in Figure 16a\cite{[M1]}.  But, this would create an extension of the
disc $\Delta_{m^o}$ that contains a closed loop in its $G_{-,-}$ graph.
Again, by Lemma 3.8 of \cite{[BF]} such loops cannot exist.  Thus,
$m^o$ can only intersect the $G_{\e,\d}$ \& $G_{-\e,-\d}$ each once.  Statement (d)
is now established and, thus, our lemma is proved.
\qed

We now finish our argument that torus knots are exchange reducible.
The main idea is simple.  Given a torus knot on $\ta$ we consider the annulus or
torus-minus-a-disc $\cS_K$.  If every possible Figure 16a\cite{[M1]} type arc $\g \subset m^i$
cannot be used to perform a type Figure 15b\cite{[M1]} alteration to the foliation of $\ta$,
it is because $\g$ essentially intersects $K$ in $\cS_K$.  So
$\g$ must locally cuts across, say, $G_{\e,\d}$ and locally $K$ is parallel to
$G_{\e,\d}$.  But, by statement (d) of Lemma \ref{lemma:the intial foliation of Delta_m^o}
$m^o$ will be locally parallel to $G_{\e,\d}$.  We can then use the edge-path in
statement (c) of Lemma \ref{lemma:the intial foliation of Delta_m^o} to produce, after
a sequence of exchange move and change of fibrations, a valence-one vertex as
depicted in Figure 16a\cite{[M1]} (except in will be in the resultant foliation of $\Delta_{m^o}$).
Since the $\g$ in $m^o$ will locally be either parallel to $K$ or not in $\cS_K$, we can
perform the alteration to the foliation of $\ta$ that is depicted in $\ta$.

%

\section{Continued analysis of standardly tiled tori and proof of Theorem
\ref{theorem:addendum theorem}}
\label{section:addendum analysis}
We now continue our analysis of the standard tiling $\ta$ of the cabling torus
that in a neighborhood of $\cS_K$ carries $K_{(P,Q)}$.  Our goal in this section is achieving a total understanding
of the embedding of such cabling torus in $S^3$.

The key to this understanding
is based upon an understanding of the meridian curves of $\ta$.
From the analysis in \cite{[M1]}, and now \S \ref{section:mistake}, we know
that all of the singularities in the foliation of $\Delta_m$ adjacent to the boundary curve
$m$ are of the same parity except for one singularity, e.g. see Figure 14a\cite{[M1]}.
It is convenient to assume the case where $\partial \cS_K$ is equal to either
$G_{+,+} \cup G_{-,-}$ (in the case where we have an annulus), or
$\{ {\rm edge-path \ of \ } G_{+,+} \} \cup \{ {\rm edge-path \ of \ } G_{-,-} \} \cup \{ {\rm edge \ of \ }G_{+,-} \} \cup \{{\rm edge \ of \ }G_{-,+} \}$
(in the case where we have a torus-minus-disc).  (Figure \ref{figure:staircase}
incorporates this assumption.)
This implies that as we can position
$m$ so that it has a {\em simply zig-zag loop}: that is, as we traverse $m$
it alternates between intersecting $G_{+,+}$ and $G_{-,-}$; and it intersects
both $G_{+,-}$ and $G_{+,-}$ once.  (Again, see Figure \ref{figure:staircase}.)
Observe that since $m$ as
simply zig-zag loop is transversely intersects the leaves of $\ta$ then $m$ is
a braid representation of the unknot in the fibration $\fibr$.

In order to help the reader visualize our ultimate goal it is useful to momentarily
put aside our foliation machinery for $\ta$ and give a description of how a cabling torus
in $S^3$ can be embedded.  (We will establish in Theorem
\ref{theorem:all cabling tori are standard} that our description here is in fact
the general case.)

\begin{figure}[htpb]
\centerline{\includegraphics[scale=.55, bb=0 0 588 377]{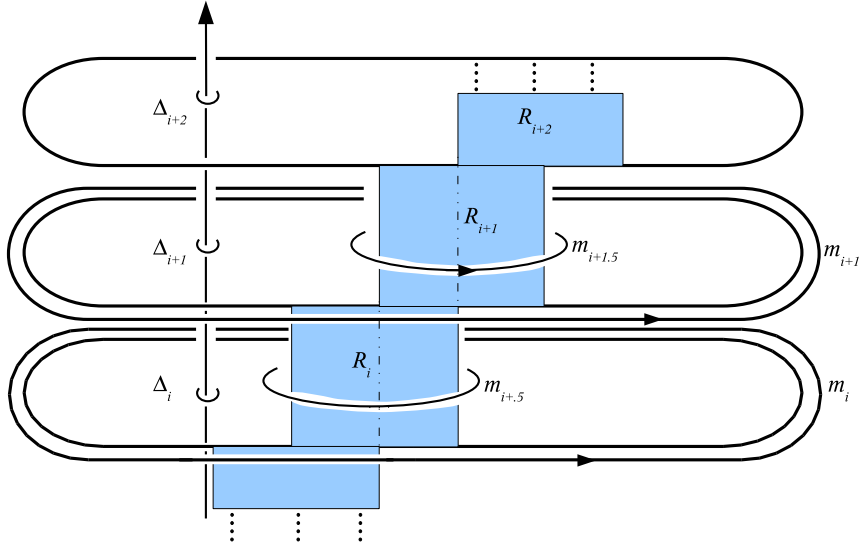}}
\caption{\footnotesize The steps configuration has the $R_i$ rectangles
arranged so that each``step'' has a $\Delta_i$ attached.
The boundary of a regular neighborhood of a steps configuration is
a standard cabling torus.  The illustration depicts curves $m_{i + .5}$,
$m_{i+1}$ and $m_{i + 1.5}$ which are meridian curves on the standard cabling torus.  (The
standard cabling torus is not explicitly depicted.)}
\label{figure:discs plus blocks}
\end{figure}

First we describe a {\em steps configuration}, $\cR_X^s$, for a collection of 
radially foliated discs and rectangular discs in $\fibr$.
Specifically, let $\{\Delta_1 , \cdots , \Delta_l\}$ be a collection
of discs that are arranged in a braid fibration setting so that each $\Delta_i$
is punctured once by the
braid axis $\axis$ and has its boundary curve transverse to the braid fibration $\fibr$. (Thus, each $\Delta_i$ is radially foliated.)  Moreover,
suppose that the indexing variable is such for any triple
$(\Delta_i, \Delta_{i+1} , \Delta_{i+2} )$
the axis $\axis$ cyclically intersects these discs in the same order.
(We treat $\{\Delta_1 , \cdots , \Delta_l\}$
as a cyclic ordering of discs.)  To continue our description of the stepping configuration, let
$\{R_1, \cdots , R_l \}$ be a collection of rectangles such that:
\begin{itemize}
\item[(i)] $\partial R_i = {r_i}^1 \cup {r_i}^2 \cup {r_i}^3 \cup {r_i}^4$.
\item[(ii)] $R_i$ is attached to $\partial \Delta_i$ along ${r_i}^2$.
\item[(iii)] $R_i$ is attached to $\partial \Delta_{i+1}$
along ${r_i}^4$.
\item[(iv)] ${r_i}^1$ \& ${r_i}^3$ are contained in disc fibers of $\fibr$
with ${r_i}^1$ the {\em bottom side} and $r_i^3$ the {\em top side} (as understood
in condition (1a) of our definition of rectangular block presentation).
\item[(v)] The leaves of the induced foliation of $R_i$ by $\fibr$ are
arcs the have their endpoints on ${r_i}^2$ and ${r_i}^4$ that are parallel
to ${r_i}^1$ and ${r_i}^3$.
\item[(vi)] The set $R_1 \cup \cdots \cup R_l$ is connected and homotopic equivalent to $S^1$.
The associated embedding of this $S^1$ in $\reals^3$ is a braid $X$.
\item[(vii)] For some fixed integer $k >0$ we have that ${r_{i+k}}^1$ and ${r_i}^3$ are
contained in the same leaf of the induced foliation of $R_1 \cup \cdots \cup R_l$.
\end{itemize}

We refer the reader to Figure \ref{figure:discs plus blocks}.  The reader should take notice
of the ``steps'' arrangement of the set $R_1 \cup \cdots \cup R_l$ and that by construction
each ``step'' of the configuration has a $\Delta_i$ disc attached.  Thus, a regular neighborhood
of $(\cup_{1\leq i \leq l} \Delta_i) \cup (\cup_{1 \leq j \leq l} \Delta_j)$ is topologically
a solid torus.  We will call the boundary of such a neighborhood a {\em standard cabling torus}.
(Referring back to our exacting description of Theorem \ref{theorem:addendum theorem} immediately after the statement
of Theorem \ref{theorem:addendum theorem}, we are working towards showing that the peripheral torus
of $K(P^\prime , Q^\prime )$ is always a standard cabling torus.)  Finally, restricting to
the collection of rectangles, the reader should notice that $\cR_X^s$ contains an underlying
rectangular block presentation $\cR_X$.

Although we will not appeal to the concept until \S\ref{section:addendum theorem}, we say that
$\cR_X^s$, and its associated standard cabling torus $\ta$, has {\em interlocking
homogeneous twisting} if the underlying rectangular block presentation $\cR_X$ does.

\begin{figure}[htpb]
\centerline{\includegraphics[scale=.6, bb=0 0 433 295]{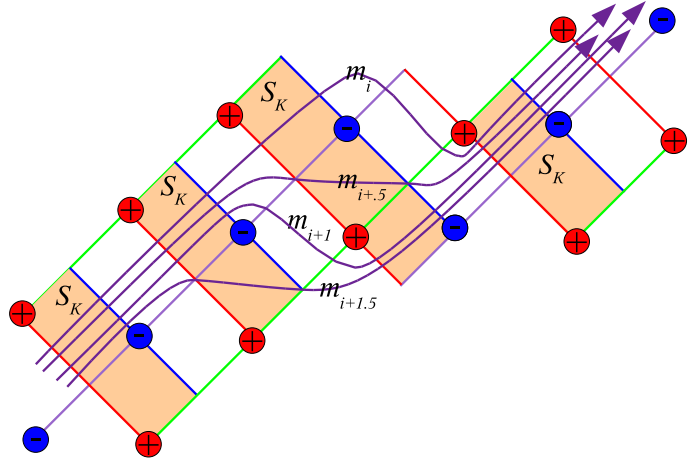}}
\caption{\footnotesize The tiling foliation of $\ta$ can be ``laid out flat'' in a
staircase pattern.  The limitation that $m$ intersect exactly once $G_{+,-}$ and
$G_{+,+}$ illustrated.  (A choice of parity is made for the purpose of illustration
simplification since there is a similar figure using $G_{-,+}$ and $G_{-,-}$.)}
\label{figure:meridian walk}
\end{figure}

It is readily seen that a standard cabling torus is in fact one of our tiled $\ta$:
each vertex is valence four and the parity of the tiles correspond to a checker-board pattern.
Specially, we refer the reader to the definition of ${\bf b}$-support of $K$, $\cS_K$
and the statement of Proposition 5.8\cite{[M1]} governing the behavior of $\cS_K$.
By Proposition 5.8\cite{[M1]} topologically $\cS_K$ is either an annulus or a
torus minus an open disc.  In either case, by observation-$\star$
we now know that a meridian curve transversely cuts through
$\cS_K$ everywhere except at one place where it must intersect, say, $G_{+,-}$ and $G_{-,+}$ once.
In Figure \ref{figure:meridian walk} we have depicted the meridian curves $m_i$ \& $m_{i+1}$
which cut across $\cS_K$ several times, transversally intersecting $G_{+,+}$ and $G_{-,-}$, but
only once are depicted as intersecting $G_{+,-}$ and $G_{-,+}$.  The labeling of the curves $m_i$ \& $m_{i+1}$ in
Figure \ref{figure:meridian walk} is meant to correspond to the labeling of curves
$m_i$ \& $m_{i+1}$ in Figure \ref{figure:discs plus blocks}.  We ask the reader to notice that in
Figure \ref{figure:discs plus blocks} all of the $m$-curves are isotopic
to each other in the complement of
the standard disc-rectangle configuration, and thus are isotopic meridian
curves on the standard cabling torus.
Moreover, we ask the reader to see that the isotopy from $m_i$ to $m_{i+1}$
which has $m_{i+ .5}$ as an
intermediary curve is in fact an exchange move of type-I
as illustrated in Figure 11\cite{[M1]}: the $n$-braid $m_i$ passes through the
axis $\axis$ to become $m_{i+.5}$, a meridian curve that has
two points of tangency with the braid fibration
$\fibr$; and the passes through $\axis$ to become the $n$-braid $m_{i+1}$.
The corresponding sequence is
depicted in Figure \ref{figure:meridian walk}: The curve $m_i$ is
transverse to the foliation of $\ta$; the
isotopy of $m_i$ to $m_{i+.5}$ passes through a negative vertex
and $m_{i+.5}$ intersects two singular points
in the foliation---points where $m_{i + .5}$ is necessarily tangent
to disc fibers of $\fibr$; and, finally
$m_{i+.5}$ passes through a position vertex to become $m_{i+1}$
which is transverse to the foliation of $\ta$.
Iterating this isotopy on the new $n$-braid $m_{i+1}$ we can swipe
out the entire foliation of $\ta$ until we arrive back at $m_i$.

\begin{figure}[htb]
\centerline{\includegraphics[scale=.65, bb=0 0 467 372]{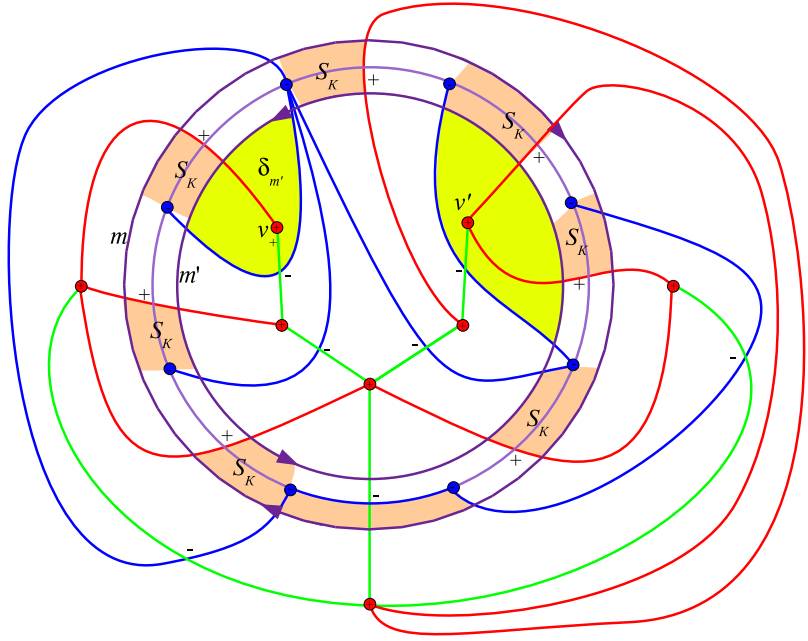}}
\caption{\footnotesize The meridian curves $m$ \& $m^\prime$ are simply zig-zag loops when
viewed in the foliation of $\ta$, and they co-bound an annulus $\cA \subset \ta$ as illustrated.
The portion of $ \cS_K \cap \cA \subset \cA$ is indicated by the shaded regions above.
Attached to $\cA$ along its boundary curves are the two meridian
discs $\Delta_m$ \& $\Delta_{m^\prime}$.  A sub-disc $\delta_{m^\prime} \subset \Delta_{m^\prime}$
has $\partial \delta_{m^\prime} = \a_1 \cup \a_2$ where
$\a_1 \subset m^\prime$ \& $\a_2 \subset G_{-,-}$; and 
$\delta_{m^\prime}$ contains exactly one vertex of $v_+$ which is
an endpoint of $G_{+,+}$.  The vertex $v^\prime$ illustrates an opportunity for a non-standard
change of fibration.}
\label{figure:sphere 2}
\end{figure}

As previously stated, our goal is to show that all the $\ta$
of any triple $(K, \ta ,m)$ is a standard cabling
torus.  In order to show this it is sufficient to show that there
is a meridian curve $m_i$ satisfying
observation-$\star$ that is a $1$-braid with respect to $\axis$.
Thus, we have the following proposition.

\begin{prop}
\label{proposition:braid index of meridians}
Let $m \subset \ta$ be a meridian curve position as a simply zig-zag loop.
Assume that $m$ intersects $2{\rm s}$ singular leaves in the foliation of $\ta$.
Then the braid index of $m$ is either $1$ or ${\rm s} - 1$.
\end{prop}

\pf
We start by considering two parallel meridian curves, $m ,m^\prime \subset \ta$, which co-bound
a sub-annulus $\cA \subset \ta$ which has an induced foliation that contains only negative
vertices of $\ta$.  (See Figure \ref{figure:sphere 2} and its caption description.)
To elaborate further,
the foliation of $\ta$ restricted to $\cA$ will have ${\rm s}$ positive vertices and ${\rm s}$
singular points.  Since $m$ \& $m^\prime$ are simply zig-zag loops
all but one of these singular points will have positive parity.

We now use that fact that $m$ and $m^\prime$ both bound meridian discs in $\ta$ to extend
$\cA$ to a $2$-sphere.  Let $\Delta_m$ and $\Delta_{m^\prime}$ be meridian disc in $\ta$ having
having $m$ and $m^\prime$ as boundary curves, respectively.  Then
$\Sigma = \Delta_m \cup \cA \cup \Delta_{m^\prime}$ is a $2$-sphere.  As previously noted, the foliations
of $\Delta_m$ and $\Delta_{m^\prime}$ will only positive vertices in their foliation, i.e. all
of the negative vertices live in $\cA$.

With this setup we start our argument with the claim that in $\Delta_{m^\prime}$ there exists
a sub-disc $\delta_{m^\prime}$ such that: $\partial \delta_{m^\prime} = \a_1 \cup \a_2$ where
$\a_1 \subset m^\prime$ \& $\a_2$ is in the $G_{-,-}$ graph of $\delta_{m^\prime}$;
and, $\delta_{m^\prime}$ contains exactly
one vertex of the $G_{+,+}$ (of $\Delta_{m^\prime}$) which is an endpoint of $G_{+,+}$.
The arc $\a_2$ is necessarily
then contained in a negative singular leaf of the foliation of $\Sigma$ and $\a_1$ must
have angular length $2\pi$ in $\fibr$.  To see this first notice that by construction
$\Sigma {\rm 's}$ $G_{-,+}$ graph is totally contained in $\cA$; and, its $G_{+,-}$ graph must then intersects
$\cA$ in exactly one arc that cuts through the single negative singularity contained in
$\cA$.  Since all of $\Sigma {\rm 's}$  $G_{\e,\d}$ graphs must be simply connected
(a property of such graphs on $2$-spheres following again from Lemma 3.8 of \cite{[BF]}),
we have $\Delta_{m^\prime}$ containing an endpoint of $G_{-,+}$ which we label $v_+$.  
If $v_+$ is the only positive vertex in $\Delta_{m^\prime}$ then from $v_+$ we can see
$m^\prime$ in sequence: intersects $G_{+,+}$; pass through $\cS_K$; intersects $G_{-,-}$;
and, finally, an intersection with $G_{+,+}$.  Thus, we can perform a non-standard change of
fibration between the middle $G_{-,-}$ singularity and the last $G_{+,+}$ singularity.
So assume that $\Delta_{m^\prime}$ contains more than one positive vertex.  Now again suppose
that $v_+$ is adjacent to more than one edge of $G_{+,+}$.  This again would imply that
from the vertex $v_+$ we can see in sequence: intersects $G_{+,+}$;
pass through $\cS_K$; intersects $G_{-,-}$; and an intersection again with $G_{+,+}$.
But, this again would give us the ability to perform a non-standard change of fibration
to alter the foliation of $\ta$ so as to produce valence three, then valence two vertices.
(Recall valence two vertices can be eliminated using exchange moves.)

\begin{figure}[htpb]
\centerline{\includegraphics[scale=.65, bb=0 0 661 312]{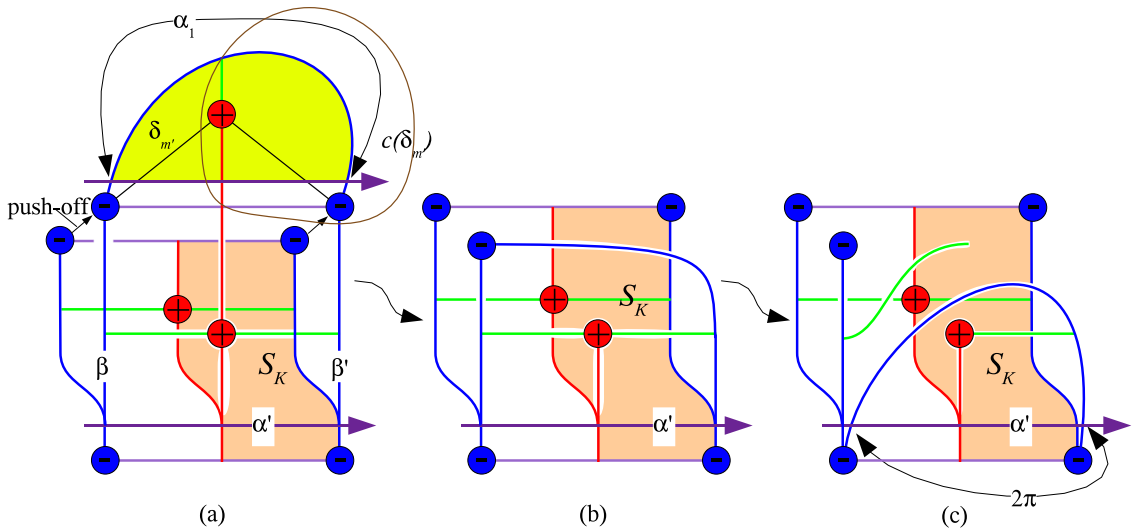}}
\caption{\footnotesize In the illustration the push-off of $\delta_{m^\prime} \cup (R \setminus \a^\prime)$
is shown ``above'' the foliation of $\ta$.  In (a) we perform an exchange move predicated on the
existence of the valence two vertex in $\delta_{m^\prime}$ so that the $\bb$-arcs inside the
$c(\delta_{m^\prime})$-curve are inessential an can be eliminated to produce the foliation in (b).
We go from (b) to (c) via a standard change of fibration.}
\label{figure:shifting a_1}
\end{figure}

So we now have an $\a_1 \subset m^{\prime}$ that has angular length $2\pi$ and, when
transversed in the direction of it orientation, intersects
$G_{+,+}$, passes through $\cS_K$, then intersects $G_{-,-}$.
(See Figure \ref{figure:shifting a_1}.)
We now refer to Figure \ref{figure:shifting a_1} and use its configuration to argue that as we move along
$\cS_K$ there are arcs having similar intersection sequence that also angular length $2\pi$.
Specifically, we focus on the sub-disc $\delta_{m^\prime}$ that is attached to $\ta$ along
the arc $\a_1$.  Let $R$ be a rectangular region in $\ta$  such that:
$\partial R = a_1 \cup \b \cup \a^\prime \cup \b^\prime$; both $\b$ \& $\b^\prime$ are arcs in the
$G_{-,-}$ graphs of $\ta$; and $\a^\prime$ is an arc whose entire is transverse to $\ta {\rm's}$
foliation and passes into $\cS_K$ after intersecting $G_{+,+}$ as shown in 
Figure \ref{figure:shifting a_1}.  We now push a copy of $ \delta_{m^\prime} \cup (R \setminus \a^\prime)$
off into the solid torus that $\ta$ bounds in $S^3$
leaving $ \delta_{m^\prime} \cup R$ attached to $\ta$ along $\a^\prime$.  (Thus, the observer's viewpoint
in Figure \ref{figure:shifting a_1} is from inside the solid torus.)  We now pass from
Figure \ref{figure:shifting a_1}(a) to  \ref{figure:shifting a_1}(b) by performing an exchange move and
an elimination of a valence two vertex.  Next, we pass from Figure \ref{figure:shifting a_1}(b) to
\ref{figure:shifting a_1}(c) by performing a standard change of fibration.  This yields an arc 
in $\ta$ that is 
slightly shifted in the foliation (which we still call $\a^\prime$), and it has the
intersection sequence of first hitting $G_{+,+}$, passing through $\cS_K$ before hitting $G_{-,-}$.
Since $\a^\prime$ intersects twice the singular leaf in what remains of the altered $R$ its
angular length must be $2\pi$.  Iterating this argument over all of $\cS_K$ we see that any oriented
arc that is: positively transverse to the foliation of $\ta$; has endpoints on $G_{-,-}$; does not intersect
$G_{-,+}$ or $G_{+,-}$; intersects $G_{+,+}$ before it enters the entire of $\cS_K$;
and has angular
length $2\pi$.  Since $m^\prime$ passes through $\cS_K$ ${\rm s} -1$-times, its angular length is at least
$2\pi ({\rm s}-1)$.  But, since the braid index of $m \cup m^\prime$ is equal to the number of vertices
contained in $\cA$, and since the braid index of $m$ is at least $1$, we must have that the
index of $m$ is $1$ and $m^\prime$ is ${\rm s}-1$.\qed

\begin{figure}[htb]
\centerline{\includegraphics[scale=.75, bb=0 0 467 372]{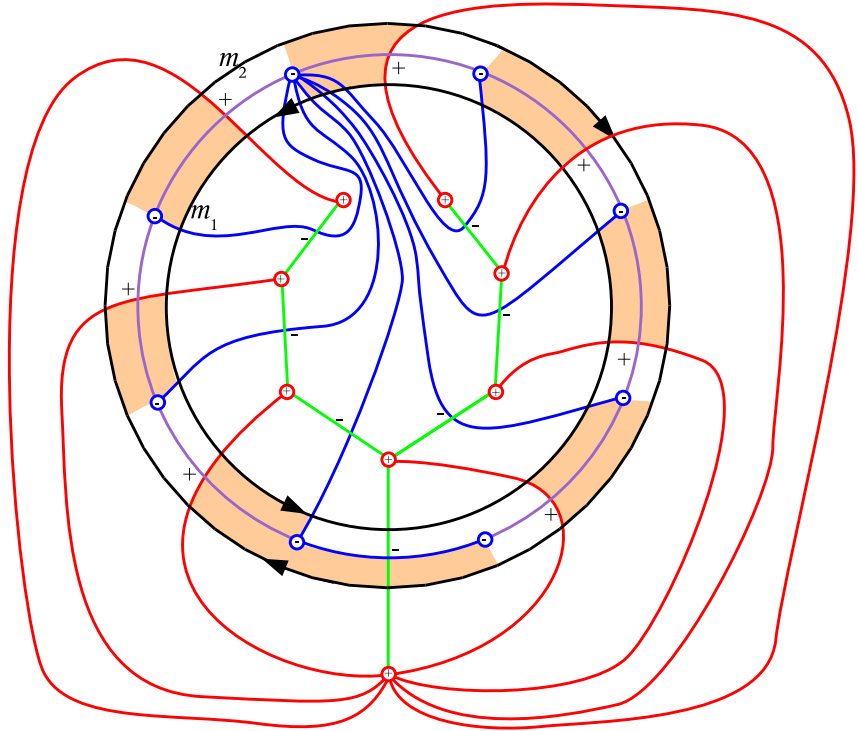}}
\caption{\footnotesize In this illustration the ``outside'' meridian disc has exactly
one vertex, thus, implying that the corresponding boundary meridian curve is a $1$-braid.}
\label{figure:actual Sigma foliation}
\end{figure}

For the propose of completeness Figure \ref{figure:actual Sigma foliation} illustrates
the foliation of a $\Sigma$-sphere that has $\Delta_m$ radially foliated and
$\Delta_{m^\prime}$ foliated with ${\rm s}-1$ positive vertices (where ${\rm s}=8$).

We are now in a position to establish our previously stated goal.

\begin{thm}
\label{theorem:all cabling tori are standard}
The $\ta$ of any triple
$(K, \ta ,m)$ is a standard cabling torus.
\end{thm}

\pf
We will establish that $\ta$ is a standard cabling torus by building the associated
$\cR_X^s$.  To do this we need to specify the $\Delta_i$-discs and the $R_j$-rectangles.

Let $\{m_1, \cdots , m_l\}$ be a sequence of zig-zag meridian
curves all of braid index one such that $m_i+1$ is obtained from $m_i$ via the type-I
exchange move illustrated in Figure \ref{figure:meridian walk}, and $m_1$ is obtained
from $m_l$ also by such a move.  Then as we push $m_1$ to $m_2$ to $m_3$ until we
arrive back at $m_1$ we will have sweep out our torus $\ta$.  Now let $\Delta_i$ be a
meridian disc inside the solid torus that $\ta$ bounds.  Since $m_i$ is a one braid
$\delta_i$ will be punctured algebraically once by $\axis$ and, thus, geometrically once
by $\axis$.  (That is, exchange moves and the elimination of any inessential $\bb$-arcs
in the foliation of $\Delta_i$ will produce a radially foliated discs.)  We let
this collection of radially foliated discs be our radially foliated discs of our
steps configuration.

\begin{figure}[htb]
\centerline{\includegraphics[scale=.6, bb=0 0 433 295]{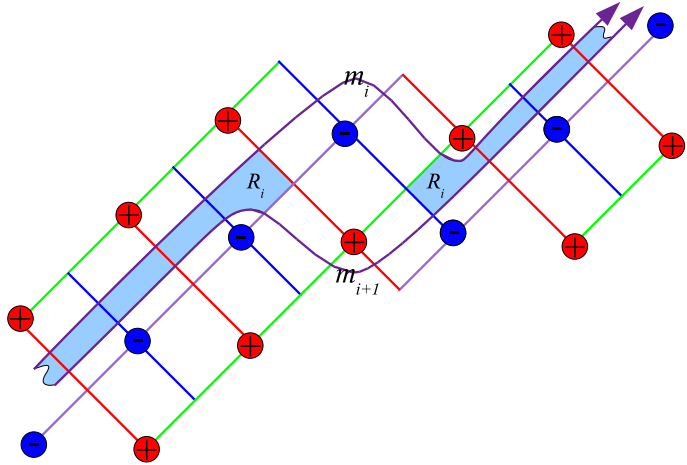}}
\caption{\footnotesize The shaded ``rectangular'' region labeled $R_i$ is attached to
$m_i$ (resp. $m_{i+1}$), and thus $\Delta_i$ (resp. $\Delta_{i+1}$).  After possibly
``shaving off'' the two sides contained in disc fibers of $\fibr$ $R_i$ will
then correspond to the $R_i$ in Figure \ref{figure:discs plus blocks}.}
\label{figure:choosing rectangles}
\end{figure}

It remains to describe how we can specify the rectangular discs of our needed steps configuration. 
But an initial collection of such rectangular discs is easily obtained from the
positioning of the $m_i {\rm 's}$ in the foliation of $\ta$.  Referring to Figure
\ref{figure:choosing rectangles} (which is a variation on Figure \ref{figure:discs plus blocks})
we see that there is a natural region between $m_i$ and $m_{i+1}$ that is between the
two singular leaves which $m_{i+.5}$ is tangent to.  This region has an induced foliation
corresponding to that of condition (v) in our description of the steps configuration.
Specifically referring to the shaded region in Figure \ref{figure:choosing rectangles}:
$r_i^2 = R_i \cap m_i$; $r_i^4 = R_i \cap m_{i+1}$; $r_i^1$ is the edge-path in one of
the singular leaf in $\partial R_i$; and $r_i^3$ is the edge-path in the other
singular leaf in $\partial R_i$.  Since singular leaves are contained in disc fibers
of $\fibr$ we have the $r_i^1$ \& $r_i^3$ are contained in disc fibers.  
By choosing $r_i^1$ to be associate with the negative singularity and $r_i^3$ associated
with the positive singularity we can have $r_i^1$ being the bottom side and $r_i^3$
being the top side of our $R_i$ rectangle.  So conditions
(i)-(iv) of our steps configuration description are satisfied.  By ``shaving off'' the
$r_i^1$ and $r_i^3$ sides of all the $R_i$ rectangular discs we can easily obtain
condition (v).  Now we notice that the union of all these $R_i {\rm 's}$ in $\ta$
is topologically an annulus so we obtain condition (vi).  Condition (vii) follows
from the regularity of the staircase pattern of our $\ta$ foliation in Figure
\ref{figure:staircase}.  
Specifically, the fixed $k$ in condition (vi) is equal to $s-2$ where $2s$ is the number
of singular leaves a zig-zag meridian curve intersects in the foliation of $\ta$.
Clearly the boundary of a regular neighborhood of
$\cR_X^s = (\cup_{1 \leq i \leq l} \Delta_i) \cup (\cup_{1 \leq j \leq l} R_j)$ is again a standardly tiled
torus that is isotopic to our original $\ta$ in our axis/fibration coordinate system.
\qed

%

\section{Proof of Theorem \ref{theorem:addendum theorem}}
\label{section:addendum theorem}

We restate the argument in \cite{[M1]} up to the error in Theorem 1.1.  We begin with a quadruple
$(K,\ta,m,\Delta_m)$.  Through a sequence of exchange moves
and destabilizations we can assume that $\ta$ (which contains $K$) has been position
arbitrarily close (in the $(\axis,\fibr)$ coordinate system) to a knot $K^\prime$
that is positioned on a cabling torus $\ta^\prime$ that has either:
a circular foliation; or, a mixed foliation; or, a tiled foliation.
For the first two possibilities by Corollary 3.1\cite{[M1]} and Proposition 4.1\cite{[M1]} we
have that $K$ is exchange reducible.  For the tiled foliation case we are left with
$K^\prime$ being positioned in the $\bb$-support $\cS_{K^\prime} \subset \ta^\prime$
that is either
an annulus or a torus-minus-a-disc, as described in Proposition 5.8\cite{[M1]}.
To establish our theorem we need to show that
the associated $\cR_X^s$ to
$\ta^\prime$ has its underlying rectangular block presentation $\cR_X$
having interlocking homogeneous twisting.
As the reader may now suspect we will derive the needed
rectangular block presentation from the tiled foliation of $\ta^\prime$.

So our starting point is in fact the tiled foliation of $\ta^\prime$.  As described in
the proof of Theorem \ref{theorem:all cabling tori are standard} from the tiled foliation
of $\ta^\prime$ we can readily produce a step configuration for which the boundary of
a regular neighborhood corresponds to $\ta^\prime$.  By positioning each
$\Delta_i$ disc in the configuration to be contained in an appropriate level
plane $P_{z_0}$ and scaling their radii to be $1$, we can achieve the condition that
all of the $r_i^2$- \& $r_i^4$-sides of the configuration's rectangular discs are contained
in $\cT_1$.  Moreover, we can then easily isotop the rectangular discs to have
empty intersection with the interior of the solid cylinder that $\cT_1$ bounds in $\reals^3$.
Once this positioning of the rectangular discs it is self evident that their union is
a rectangular block presentation $\cR_X$ of the braid $X$ which passes through each
$R_i (\subset \ta^\prime)$ once and, thus, intersects each simply zig-zag meridian curve on
$\ta^\prime$ once.  Clearly $X$ is a longitude of $\ta^\prime$;
$\cup_{1\leq i \leq l} R_i = \cR_X$ is a rectangular block presentation; and $K^\prime$ is
a cabling of $X$.

We now need to establish that $\cR_X$ has homogeneous twisting and is interlocking.  The argument
for each is to show that if $\cR_X$ does not possess the claimed property then there will
exist a standard change of fibration which is not obstructed by $K^\prime$ that allows us
to continue to simplify the foliation of $\ta^\prime$ through braid isotopies, exchange
moves and $\pm$-destabilizations of $K^\prime$.  
We use two claims to organize the remaining argument.

\begin{figure}[htpb]
\centerline{\includegraphics[scale=.6, bb=0 0 490 153]{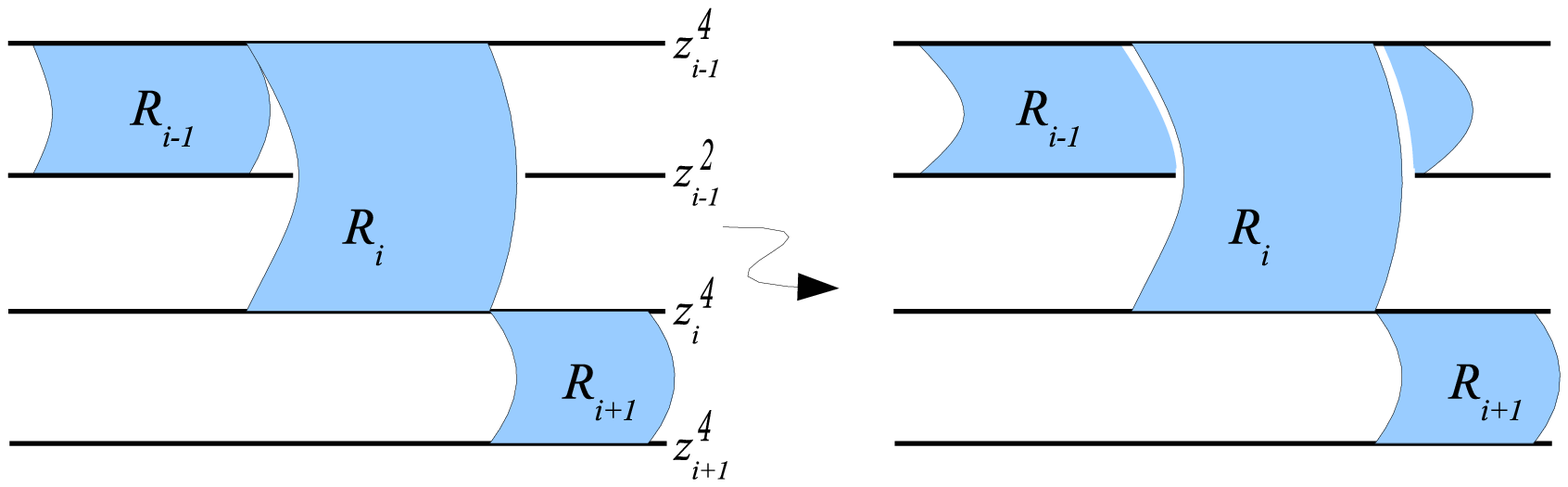}}
\caption{\footnotesize }
\label{figure:block slide}
\end{figure}

\noindent
\underline{Claim 1}--$\cR_X$ has homogeneous twisting.

We refer the reader to Figure \ref{figure:block slide} for this argument.
Reiterating the condition for homogeneous twisting, we require that for any
triple $(z_{i-1}^4, z_i^4, z_{i+1}^2)$ satisfying the condition that
$z_{i+1}^4 < z_{i-1}^4 < z_i^4$ then we have that $ z_{i+1}^4 < z_{i+2}^4 < z_i^4$.
Suppose this is not the case. Then either $ z_{i+2}^4 < z_{i+1}^4 < z_i^4$ or
$ z_{i+1}^4 < z_i^4 < z_{i+2}^4 $.  The left-hand side illustration
in Figure \ref{figure:block slide} depicts
the situation when $ z_{i+1}^4 < z_{i+2}^4 < z_i^4$.  We are allowing this figure
to serve a dual purpose:  it illustrates a portion of a rectangular block presentation $\cR_X$;
and it illustrates a portion of a step configuration $\cR_X^s$.  When viewing
a step configuration one should realize that
associated with each block is a negative singularity (before the $r_j^1$-side in the
fibration $\fibr$) and a positive singularity (after the $r_j^3$-side in $\fibr$).
Now going back to our current figure,
if we can ``slide'' the $r_{i-1}^3$-side of $R_{i-1}$ forward in $\fibr$ along
$\partial \Delta_{i-1}$ and $\partial \Delta_i$ so as to move it past the $r_i^1$-side
of $R_i$ then in the corresponding tiling of $\ta^\prime$ we will have preformed a
nonstandard change of fibration to the tiling.  One should observe that the
braid $K^\prime \subset \ta^\prime$ will not act as an obstruction to such a
nonstandard change in fibration since $K^\prime$ can be viewed as being positioned within
an arbitrarily small neighborhood of the set 

\begin{equation}
(\cup_{1\leq i \leq l} (\partial \Delta_i))
\cup (\cup_{1 \leq i \leq l} \l_i) \subset (\cup_{1 \leq i \leq l} \Delta_i) \cup (\cup_{1\leq i \leq l} R_i) \subset \ta^\prime .
\label{equation:Kprime set}
\end{equation}
So the only possible obstruction to performing this slide will come from the positioning
of the rectangular blocks themselves.  However, as seen in Figure \ref{figure:block slide},
the $R_i$ block ``shields'' the sliding of the $r_{i-1}^3$-side of the $R_{i-1}$ block
from obstruction---all of the rectangular blocks positioned ``underneath'' $R_i$ can
be pushed forward in $\fibr$ before we slide $r_{i-1}^3$ past $r_i^1$.
Thus, our claim is established.  (A similar illustration
can be used to depict the situation where $ z_{i+1}^4 < z_i^4 < z_{i+2}^4 $.  This
argument can also be presented in a $H_\theta$-sequence form.)

\noindent
\underline{Claim 2}--$\cR_X$ is interlocking.

The argument establishing this claim is almost a tautology.  As set equation
\ref{equation:Kprime set} specifies, whatever alteration we make to a step configuration
we cannot disturb the $(\cup_{1\leq i \leq l} (\partial \Delta_i))
\cup (\cup_{1 \leq i \leq l} \l_i)$ set.  So if we try to forward slide the $r_j^3$-side of the
$R_j$ rectangular block along $\partial \Delta_j \cup \partial \Delta_{j+1}$ past the
$r_{j+1}^1$-side of the $R_{j+1}$ block the only possible obstruction will be some
$r_{u_j^1}$-side of a rectangular contained in a
$\l_{u_j^1} (\subset \cup_{1 \leq i \leq l} R_i)$ leaf.  So we can try to slide
$r_{u_j^1}^3$ along with all of $\l_{u_j^1}$ forward and past $r_j^1$.  But, again
we made be obstructed by another $r_{u_j^2}^3$ that is contained in a
$\l_{u_j^2}$ leaf.  Continuing in this manner either we can push all obstructions to
pushing $r_j^3$ past $r_{j+1}^3$ or we find that $r_{j+1}^3 (\subset \l_{j+1}^3)$ itself
is an obstruction to pushing $r_j^3$ past $r_{j+1}^3$.  It is the latter that gives us
an interlocking rectangular block presentation.  If we can push all of the obstruction
out of the way then we can again perform a non-standard change of fibration to
the tiled foliation of $\ta^\prime$ and reduce the complexity of its foliation
through a sequence of braid isotopies, exchange moves and $\pm$-destabilizations.
this establishes our remaining claim and ends the proof of Theorem
\ref{theorem:addendum theorem}.\qed

%

\section{Building interlocking homogeneous twisting rectangular block presentations.}
\label{section: building intlocking presentations}

In this section we will give a general procedure for building the interlocking homogeneous
twisting step configuration $\cR_{K(P',Q')}^s$ on which
the knot $K_{(P,Q)}$ lives.
Necessarily we will have
$(P^\prime,Q^\prime) = \{(p_1,q_1), \cdots , (p_{h - 1},q_{h -1})\}$ for
$(P,Q) = \{(p_1,q_1), \cdots , (p_{h},q_{h})\}$.
(Our notation $(P,Q)$ and $(P^\prime, Q^\prime)$
is consistent with that in \S\ref{section: Introduction}.)

This procedure is described in terms of the
$H_\theta$-sequence and can then be readily
geometrically realized as a the step configuration.
Next, using the geometrically realized of the step configuration we can
superimpose over it a {\em rectangular diagram} of a
{\em rectangular diagram} of the cable knot $K_{(P,Q)}$.  Specifically,
$K_{(P.Q)}$ will be comprised as an union of {\em horizontal arcs} contained in the unit cylinder
$\cT_1$; and {\em vertical arcs} contained in half-plane $H_\theta \in \fibr$ and having
$\rho \geq 1$.  

There is an immediate advantage to presenting $K_{(P,Q)}$ as a rectangular diagram.
A rectangular diagram readily yields a transversal braid presentation and a Legendrian
knot presentation of the underlying knot type in the standard symmetric contact structure
of $\reals^3$, i.e. the plane field coming from the kernel of the $1$-form
$\alpha = d z + \rho^2 d \theta$.

In general, by slightly tilting vertical arcs and smoothing corners any rectangular diagram
of a knot $K$
can be deformed into a closed braid.  Bennequin's classical
result \cite{[B]}, that a closed braid can also be thought of a a transversal link in the standard
contact structure of $\reals^3$ or $S^3$, gives us that after this deformation we are looking
at a transversal presentation of $K$.  The classical invariant of transversal knots,
the Bennequin number, is readily computed by taking the algebraic
length less the braid index.  That is, $\beta ( K) = \ell - n $.

Seeing how rectangular diagrams of a knot $K$ can be thought of as Legendrian presentations is almost as
straight forward.  In general, a rectangular diagram can initially be associated with a
piecewise Legendrian diagram:  The horizontal arcs can be associated with leave segments in
the characteristic foliation of a cylinder
$C_\epsilon = \{(\rho, \theta , z) | \ \rho = \epsilon ( < 1)\}$
where $\epsilon$ can be arbitrarily small; the vertical arcs can be associated with leave
segments in the characteristic foliation of a cylinder
$C_R = \{ (\rho, \theta ,z) | \ \rho = R ( > 1) \} $
where $R$ can be arbitrarily large; and, a horizontal arc is connected to a vertical arc
by Legendrian segment that correspond to
$\{ (\rho, \theta_0 , z_0) | \  \epsilon \leq \rho \leq R \}$
where $z = z_0$ is the plane containing the endpoint of the horizontal arc and $\theta = \theta_0$
is the half-plane containing the endpoint of the vertical arc.  (So these connecting
Legendrian segments are ``suspensions'' of the corners of the rectangular diagram.)  Now smooth the corners
of this piecewise Legendrian knot to obtain a smooth Legendrian knot.  (Readers
familiar with classical presentations of Legendrian knots should realize that what we obtain
is a ``braided front projection'' of a Legendrian knot.)  The classical
invariants---the Thurston-Bennequin number $tb$ and the rotation number $r$---are
easily computed from the combinatorial information in such a diagram \cite{[Ma]}.
The Thurston-Bennequin number is the braid index less the number of upward oriented
vertical arcs---$tb(K) =  {\rm writhe}( K ) - \# ( \uparrow )$.  The rotation number is the writhe less the number of upward
oriented vertical arcs---$r ( K ) = n - \# ( \uparrow )$.

\begin{figure}[htpb]
\centerline{\includegraphics[scale=.37, bb=0 0 1043 705]{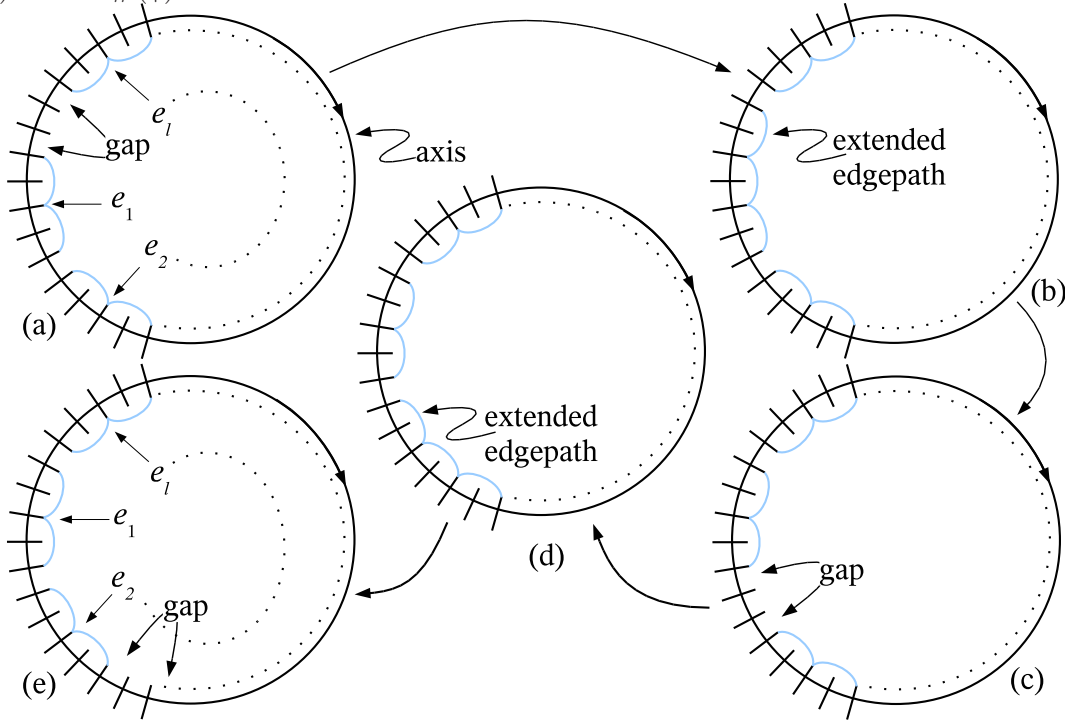}}
\caption{\footnotesize The extended edge-paths in (b) and (d) correspond to
the $\l_j$-leaves mentioned in the description of interlocking in \S\ref{section: Introduction}.
Notice that the gap is consistently ``rotating'' in the same direction.  This implies
homogeneous twisting.  Notice that the extended edge-path in (b) must occur before the
extended edge-path in (d) can occur.  This implies interlocking.}
\label{figure:discs-plus-blocks-sequence}
\end{figure}

We start with our cabling coefficient $ \{(p_1,q_1),(p_2,q_2),\cdots,(p_{h-1},q_{h-1})\}$ and
a fixed number $k$ coming from condition (3) of the description of rectangular block
presentation.  From these values we will construct an initial intersection set
$\cR_{K(P^\prime,Q^\prime)}^s \cap H_{\theta_0} \subset H_{\theta_0}$.  The reader should consult
Figure \ref{figure:discs-plus-blocks-sequence}(a) as our construction advances.

There will necessarily be
$ [ (2k + 3)(\prod_{1 \leq i \leq h-1} p_i ) + 2 ]$ points in the set
$ \cR_{(P^\prime,Q^\prime)}^s \cap \axis$.
So we realize these points of intersection in Figure \ref{figure:discs-plus-blocks-sequence}(a)
by placing down $ [ (2k + 3)(\prod_{1 \leq i \leq h-1} p_i ) + 2 ]$ arcs that transversely
intersect the axis $\axis$.  Each such arc corresponds the two radial leaves of a $\Delta_i$-disc
used in the step configuration adjoined together at their common endpoint in $\axis$.
In Figure \ref{figure:discs-plus-blocks-sequence}(a) we will use the disc that $\axis$
visually bounds in the page as own initial $H_{\theta_0}$.  Own construction builds on
the radial leaves coming from the $\Delta_i \cap H_{\theta_0}$ that are inside this
visual disc.  The endpoints of these radial arcs inside $H_{\theta_0}$ will necessarily
be the points $\varrho -{\rm points}$ mentioned in the description of
an interlocking $\cR_X$.

Next, picking an initial
point of $z_1 \in \Delta_1 \cap \axis \subset \cR_{K(P,Q)} \cap \axis \subset \axis$
we cyclically enumerate all
there other $z_j$-points in this set as we transverse $\axis$ in the direction opposite to its orientation.  (To obtain the homogeneous feature for
$\cR_{K(P^\prime,Q^\prime)}^s$ we will need
to pay attention to the orientation direction of our labeling.  One can obtain
an alternate construction by reversing all orientations.)
We place an edge-path of length $k$ in $H_{\theta_0}$ that connects the leaves of
$(\Delta_1 \cap H_{\theta_0}), (\Delta_3 \cap H_{\theta_0}), \cdots , (\Delta_{\rm odd} \cap H_{\theta_0}),
\cdots , (\Delta_{3+2k} \cap H_{\theta_0})$.
(This edge-path corresponds to the intersection of $H_{\theta_0}$ with $k$ consecutive $R_i$
rectangles in the step configuration.)
More generally, the pattern is to connect every-other leaf in the first
$2k + 3$ leaves together by an edge-path.  So starting at the leaf associated with
$z_{2k +4}$ we again connect every-other leaf by an edge-path of length $k$.  We reiterate this
edge-path connecting pattern around the radial leaves in $H_{\theta_0}$ until we have placed
$(\prod_{1 \leq i \leq h-1} p_i)$-connecting edge-paths of length $k$ in $H_{\theta_0}$.  Thus, our
every-other pattern of length $k$ edge-paths will have accounted
for $ (2k + 3)(\prod_{1 \leq i \leq h-1} p_i )$
radial leaves in $H_{\theta_0}$ and each edge-path is a leaf in the foliation of $k$
consecutive $R_i$ rectangles in the step configuration.  Our description
of the initial $H_{\theta_0}$ is complete.

Figure \ref{figure:discs-plus-blocks-sequence} illustrates the situation when $k=1$.
The edge-paths are labeled $\{ e_1 , \cdots , e_l \}$ where $l = (\prod_{1 \leq i \leq h-1} p_i)$.
In order to obtain our homogeneous feature the direction of our labeling is again
opposite the orientation of $\axis$. 
The two-arc gap between edge-path $e_1$ and $e_l$ is also indicated.

We are now in a position to describe how the $H_{\theta}$-sequence progress so as to
allow for interlocking homogeneous twisting cabling.
The reader should consult the sequence of illustrations in
Figure \ref{figure:discs-plus-blocks-sequence}.

As we push $H_{\theta_0}$ forward in the fibration $\fibr$ we will encounter the first
leaf in the foliation of $\cR_{K(P^\prime,Q^\prime)^s}$
that intersects $k + 1$ of the $R_i$ rectangles,
i.e. where condition (vii) of the step configuration is occurs and this
is necessarily a $\l_j$-leaf coming from the description of interlocking.  This can only happen
when one of the edge-paths of length $k$ in $H_{\theta_0}$ is extended to an
edge-path of length $k+1$  (which is
a $\l_j$-leaf) by adjoining to it one of the two-arcs associated with the
gap.  After passing through the fiber containing this extended Agathe immediately
revert back to having all edge-paths being of length $k$.  However, such a sequence
must also yield homogeneous twisting while maintaining the interlock nature of
$\cR_{K(P^\prime,Q^\prime)}^s$.

To start this sequencing we focus on achieving the cabling associated with the
first pair, $(p_1 , q_1)$ of $(P^\prime,Q^\prime)$.  Specifically, we have the edge-path
$e_1$ extend to a length $k+1$ edge-path that adjoins to a radial arc associated with the
gap so that the every-other radial arc pattern is maintained.  (As mentioned, this extended Agathe
corresponds to, say, $\l_{j_1}$ leaf in the description for interlocking.)
We then have $e_1$ immediately revert back to an edge-path of length $k$.  We can do this
sequence in such a fashion that the gap is ``rotated'' past $e_1$.

Referring to the change from
Figure \ref{figure:discs-plus-blocks-sequence}(a) to \ref{figure:discs-plus-blocks-sequence}(b)
to \ref{figure:discs-plus-blocks-sequence}(c)
we see that $e_1{\rm 's}$ length goes from $2$ to $3$ back to $2$.  Notice that the gap has been
moved so as to be between $e_1$ and $e_{2}$ where its initial position was between
$e_1$ and $e_l$.  This sequence constitutes a partial positive twisting.  Continuing
we can move the gap past all edge-paths $e_{2}$ through $e_l$.
In doing so we pass through edge-paths $\l_{j_2}$ through $\l_{j_l}$.
We can then reiterate this gap moving
sequence so as to achieve twisting corresponding to $(p_1 ,q_1)$.  (We are not addressing the question of which $(p_1,q_1)$ is realizable---the
answer to this question will involve an interplay between $(P,Q)$ and $k$.)

By always rotating the two-arc gap in the same direction it is easy to check that
we have homogeneous twisting.  Simply designate a point on $\axis$ as $\{ \infty \}$ and
verify that the triple inequality for any triple $(z_{i-1}^4, z_i^4, z_{i+1}^4)$ is satisfied.
Seeing that we also have interlocking is also straight forward.  Just notice that
the edge-path $\l_{j_1}$ must disappear before the edge-path $\l_{j_2}$ can occur as
we advance through the $H_\theta$-sequence.

Now we achieve the cabling associated with the second pair, $(p_2 , q_2)$ of $(P,Q)$.
For the purpose of simplifying our indexing notation, after finishing the cabling/twisting
associate with $(p_1 , q_1)$, we relabel our edge-paths so that again the two-arc
gap is between $e_1$ and $e_l$.  Again we require that our direction of
labeling be opposite the orientation of $\axis$.
We now draw an imaginary proper arc across $H_{\theta_0}$ that is away from our
radial arcs and edge-paths which
splits off the edge-paths $e_1$ through $e_{(p_2 \times \cdots  \times p_{h-1})}$ along with the
two-arc gap.  Call this split off sub-fiber disc
$H_{\theta_0}^{2} \subset H_{\theta_0}$.  Notice that
this sub-fiber is configured similar to our initial $H_{\theta_0}$ except with
$\prod_{2 \leq i \leq l} p_i$ edge-paths (along with the two-arc gap).  We can then
repeat our twisting $H_{\theta}$-sequence only restricted to $H_{\theta_0}^{2}$.
(Again, we are not addressing the question of which $(p_2,q_2)$ is realizable.)

\begin{figure}[htpb]
\centerline{\includegraphics[scale=.75, bb=0 0 592 298]{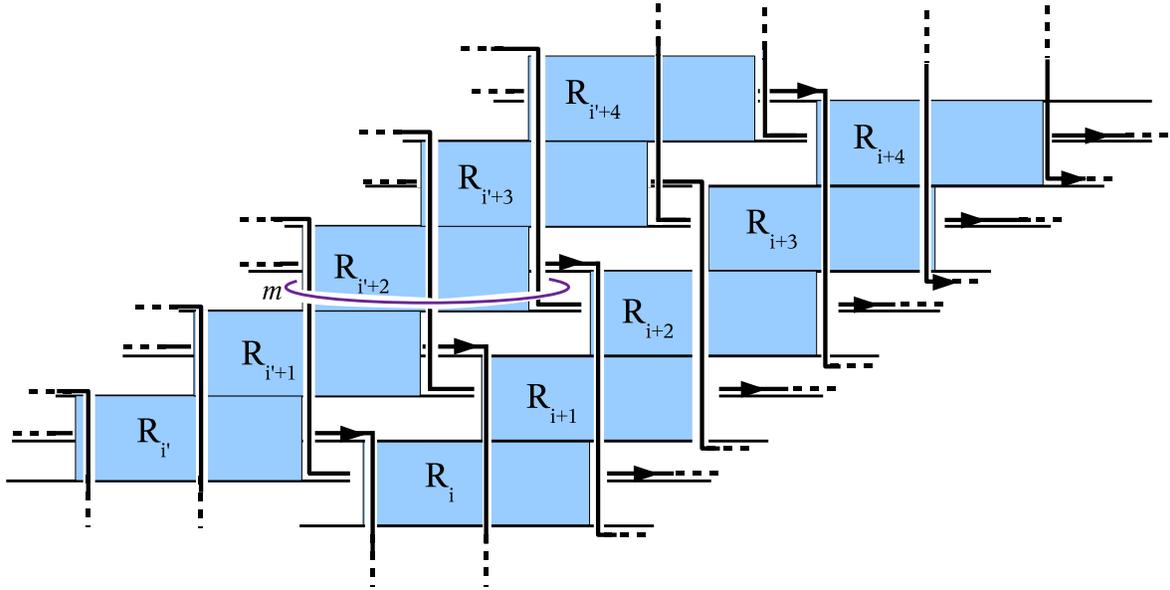}}
\caption{\footnotesize We require that within a neighborhood of every $\l_j$-leaf
there we a vertical arc of the rectangular diagram of $K_{(P,Q)}$; and, within
a neighborhood of every segment of $\partial \Delta_i$ where an $R_i$-rectangle
is not attached there be a horizontal arc of the rectangular diagram.  The meridian
curve $m$ will necessarily link $K_{(P,Q)}$ at least $k+1$ times.  Thus, $p_h \geq k+2$.}
\label{figure: imposed rectangle diagram}
\end{figure}

In general, through a sequence of sub-fiber discs---splitting off $p_3 \times \cdots \times p_l$
edge-paths, then $p_4 \times \cdots \times p_l$ edge-paths and so on---we can iterate the
underlining cabling/twisting of Figure \ref{figure:discs-plus-blocks-sequence} to obtain
the $H_{\theta}$-sequence for $\cR_{K(P,Q)}^s$.
Always relabeling the edge-paths in the opposite direction of $\axis{\rm 's}$
orientation will yield our
homogeneous twisting feature.  Again, this can be directly verified by
checking the triple inequality condition for a $z_i$-triple.
(Reversing all orientations will yield an alternative
construction.)  The every-other pattern of our edge-paths will
yield our interlocking feature because, again, the extended edge-paths correspond to
our $\l_j$-leaves.
Although care is required, a geometric
realization of $\cR_{K(P,Q)}^s$ coming from such an $H_{\theta}$-sequence can readily be produce.  A useful
format for visualizing such a realization, as we shall shown in
\S\ref{section: The Etnyre-Honda cable.}, is a square diagram that identifies the left \& right
sides and the top \& bottom sides.  (Jump ahead and see Figure \ref{figure:trefoil torus}.)

We now appeal to Figure \ref{figure: imposed rectangle diagram} to illustrate
how the knot $K_{(P,Q)}$ can be superimposed over the a geometric realization
of $\cR_{K(P^\prime,Q^\prime)}^s$.  Specifically, $K_{(P,Q)}$ will be represented utilizing a
rectangular diagram where the horizontal arcs will be line segments that are
parallel and in neighborhoods of the
$\partial \Delta_i \subset \cR_{K(P^\prime,Q^\prime)}^s$ and the vertical
arcs are in neighborhoods of the edge-path/leaves that are of length $k+1$.
For the interlocking feature to be obtain we require that there be at least vertical
arc associated with every such $k+1$ length $\l_j$-edge-path and every segment
of $\partial \Delta_i$ not having an $R_i$ rectangle attached have an associate
horizontal arc.  In Figure \ref{figure: imposed rectangle diagram} also illustrates
that $p_h \geq k+2$

In the next section we work out an explicit example.

%

\noindent
\section{Appendix: The $(2,3)$ cabling of the $(2,3)$ torus knot.  (Joint with H. Matsuda.)}
\label{section: The Etnyre-Honda cable.}

With the general construction of $\cR_{K(P,Q)}^s$ in mind along with the
corresponding rectangular diagram representation of $K(P,Q)$
we will now discuss the explicit representations to the $(2,3)$ cabling of the
$(2,3)$ torus knot (the positive trefoil) that was implicitly discovered in
Theorem 1.7 of \cite{[EH]}.  Specifically, their theorem says that there exists two
transversal knots those knot type is
$K_{((2,3),(2,3))}$ having
Bennequin number equal to $3$ but are not transversal isotopic.
The examples in this section are collaborative with H. Matsuda.

\begin{figure}[htpb]
\centerline{\includegraphics[scale=.6, bb=0 0 446 367]{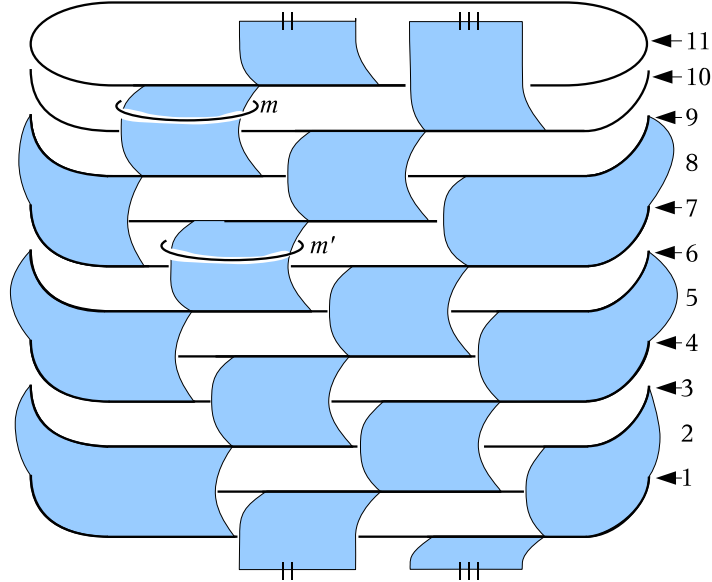}}
\caption{\footnotesize The steps configure illustrated here has three of its rectangle
discs wrap around behind the illustration and two of its rectangle discs go up
``through infinity'' and return from below.  (Or equivalently, the left and right
sides are identified; and the top and bottom are identified.)}
\label{figure:trefoil torus}
\end{figure}

From Theorem \ref{theorem:addendum theorem} we know that likely candidates for
explicit presentations could come from rectangular block presentations of the positive trefoil that
have homogeneous twisting and are interlocking.  Such a presentation was first discovered in
\cite{[P]} and is illustrated in Figure \ref{figure:trefoil torus}.  A brief pictorial narrative
is useful.  Figure \ref{figure:trefoil torus} can best be thought of initially as
the projection of $\cR_{K_{(2,3)}}^s$ onto the unit cylinder $\cT_1$ which is then
made into a torus by identifying the two ends at infinity.  Thus, the two ``half-blocks'' at the
top are identified with the two ``half-blocks'' at the bottom.  Once this understanding
is focused in one's mind the radial discs $\{\Delta_1 , \cdots , \Delta_{11} \}$ are
added to form the steps configuration.  The labels for these radial disc are the right side
column of numbers.  In addition we indicate two meridian curve, $m$ \& $m^\prime$ that
would correspond in type to the $m_{i+.5}$ meridian curves in Figure
\ref{figure:meridian walk}, i.e. in the foliation of $\ta$ they pass through a negative
and positive singularity.  Any $(2,3)$ cabling will then have to link these two
meridians three times.

Finally, referring back to conditions (3) of our definition of rectangular block presentations
and (vii) of our step configuration, we note the $k=0$ for Figure \ref{figure:trefoil torus}.
Thus, each $\l$-leaf in $\cR_{K_{((2,3),(2,3))}}^s$ is just $r_i^3 \cup r_{i+1}^1$.

\begin{figure}[htpb]
\centerline{\includegraphics[scale=.6, bb=0 0 438 357]{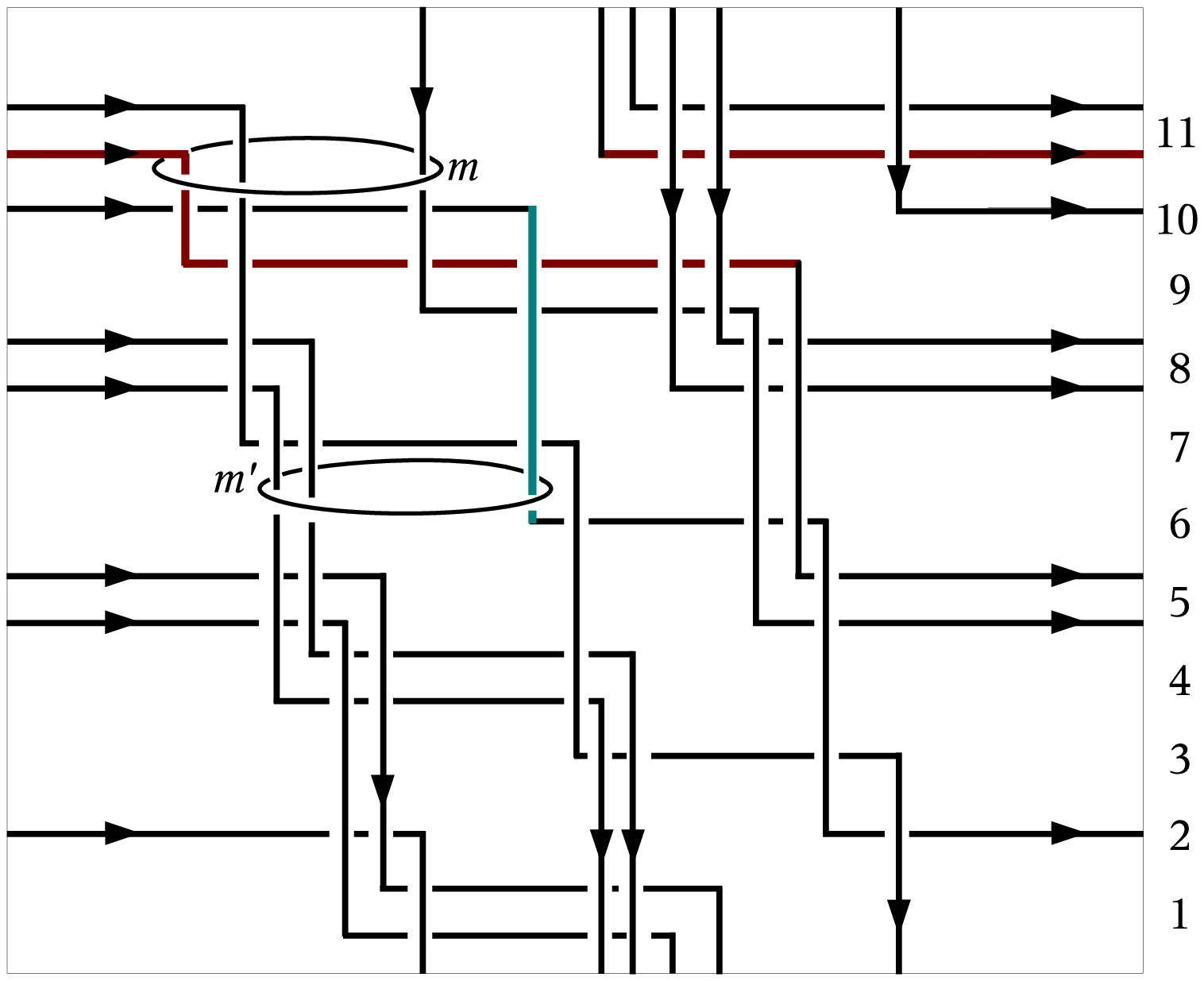}}
\caption{\footnotesize 
The illustration should be seen as a rectangular diagram projected
onto a torus.  The top side of the diagram is identified with the bottom side; and the
left side is identified with the right side.}
\label{figure:trefoil cabling}
\end{figure}

Next, Figure \ref{figure:trefoil cabling} illustrates the $(2,3)$ cabling of the torus
that is the boundary of a regular boundary of the step configuration in Figure
\ref{figure:trefoil torus}.  Our cabling $K_{((2,3),(2,3))}$ is drawn as the
rectangular diagram projection
onto $\cT_1$ with ends identified again.
The columns of numbers that
are to the right of the projection are meant to correspond to the radial disc labels
in Figure \ref{figure:trefoil torus}.  Therefore, although it makes for a very busy illustration
one could superimpose Figure \ref{figure:trefoil cabling} onto Figure \ref{figure:trefoil torus}
to understand how $K_{((2,3),(2,3))}$ is positioned in a neighborhood of $\cR_{K_{(2,3)}}^s$
and, thus, $\ta$
$\ta$.  Amongst this ambitious
visualization we have depicted the corresponding meridian curves $m$ \& $m^\prime$ in
Figure \ref{figure:trefoil cabling} that each link $K_{((2,3),(2,3))}$ three times.

We can now consider the two alternate interpretations of our rectangular diagram in
the standard contact structure.

From the Legendrian perspective we see that
Figure \ref{figure:trefoil cabling} depicts a Legendrian knot corresponding to
one describe in Theorem 1.9 of \cite{[EH]}.  Its Thurston-Bennequin number is $5$ and its
rotation number is $2$.  (See Figure 1 of \cite{[EH]}.)

If we take the associated transversal braid from Figure \ref{figure:trefoil cabling}
the Bennequin number for $K_{((2,3),(2,3))}$ is readily computed: the braid index is $8$ and
the algebraic length is $11$. Then $K_{((2,3),(2,3))}$
has a transversal knot presentation with Bennequin number $3$.

Finally, in Figure \ref{figure:trefoil cabling} we have colored two portions of the diagram.
The portion is a single vertical turquoise arc and the other is a red edge-path.  The
red edge-path has angular length in the braid a little greater than $2\pi$.  Together
these two portions indicate the existence of a negative elementary flype.  A useful
source for thinking about negative elementary flypes is \cite{[BM1]} and, in the
situation for transversal knots, \cite{[BM2]}.  Specifically, in \cite{[BM2]} it
was established that transversal knots that are presented as closed $3$-braids
and are related to each other via negative elementary flypes have the same
Bennequin number but are not transversally isotopic in the standard contact structure.
Our goal is to illustrate that the two differing transversal classes of
$TK_{(2,3),(2,3))}$ are also related by a negative flype.

\begin{figure}[htpb]
\centerline{\includegraphics[scale=.6, bb=0 0 438 357]{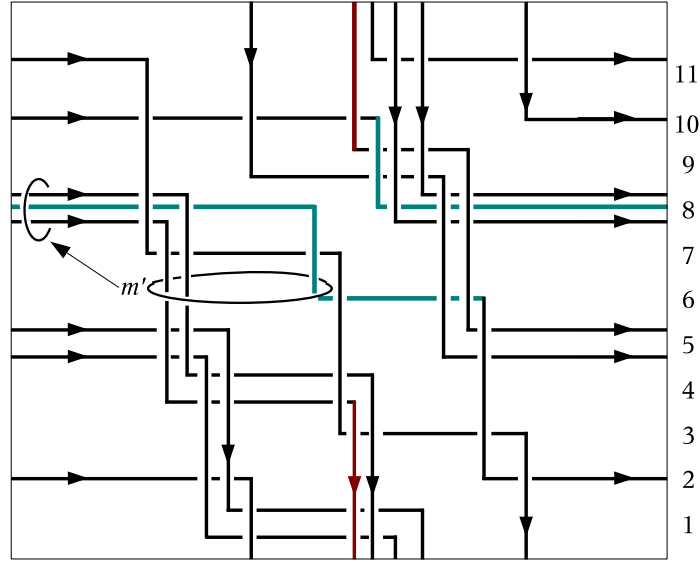}}
\caption{\footnotesize As in Figure \ref{figure:trefoil cabling}
the illustration should be seen as a rectangular diagram projected
onto a torus.}
\label{figure:trefoil cabling flype}
\end{figure}

This goal is achieved by simply depicting the rectangular diagram resulting from applying
the indicated negative elementary flype.  (See Figure \ref{figure:isotopies}(c).)
Figure \ref{figure:trefoil cabling flype} illustrates
the resulting diagram.  The turquoise edge-path corresponds to a negative stabilization
of the vertical \textcolor{Turquoise}{turquoise arc} in Figure \ref{figure:trefoil cabling};
and the \textcolor{Red}{red
vertical arc} corresponds to a negative destabilization of the red edge-path in
Figure \ref{figure:trefoil cabling}.  The meridian curve $m^\prime$ has been undisturbed by
the flype but it can now be isotopied so that it is vertical, i.e. there is some
$H_\theta \in \fibr$ such that $m^\prime \subset H_\theta$.  This implies that the resulting
cabling torus has a circular leaf and must then have a mixed foliation.  But once we
have the cabling torus possessing a mixed foliation, by Proposition 4.1\cite{[M1]}
through a sequence of braid isotopies, exchange moves and destabilizations we can produce
a foliation with only $\bc$-circles.  Thus, after the negative flype we can reduce our
cabling down to the obvious braid representation.

To connect the dots in our discussion, when viewed as closed transverse braids in the standard
contact structure, the pair of braids in Figures \ref{figure:trefoil cabling} and
\ref{figure:trefoil cabling flype} must correspond to the implicit pair in Theorem 1.7
of \cite{[EH]}.  In particular, by Theorem \ref{theorem:addendum theorem} we know that
if there is more than the transversal class that corresponds to the
exchange reducible class then such a class when viewed as a classical braid will
yield an interlocking homogeneous twisting rectangular block presentation for which it is
a cabling.  Since Figure \ref{figure:trefoil torus} illustrates such a presentation and
since by Theorem 1.7 \cite{[EH]} there are exactly two transverse isotopy classes in
$TK_{((2,3),(2,3))}$ having Thurston-Bennequin invariant $3$.  We know
Figure \ref{figure:trefoil cabling flype} must be one (the exchange reducible class)
and Figure \ref{figure:trefoil cabling} must be the other.

Finally, if we decide to view the rectangular diagram as representing a Legendrian
diagram one can readily check that its classical invariants have not be altered
(again, see \cite{[Ma]}).  So
the Legendrian ``diagrams'' of Figures \ref{figure:trefoil cabling} 
and \ref{figure:trefoil cabling flype} represent one of the points
of multiplicities in Figure 1 of \cite{[EH]}, i.e. two differing Legendrian classes that have the
same knot type and classical invariants $tb =5$ and $r=2$ as corresponding to
the multiplicity $2$ entry in
Figure 1 of \cite{[EH]}.

Three concluding remarks are in order.
First, it is interesting to note that negative elementary flypes (Figure \ref{figure:isotopies}(c))
are to date the only
method for producing explicit examples of transverse knots that are of the same knot
type having the same classical invariant, but are not transversally isotopic.
(In \cite{[BM2]} it was established that $3$-braids related by negative
elementary flypes produced such examples.)
Second, the question arises whether the MTWS braid calculus \cite{[BM1]} for cable knots
is only destabilizations, exchange moves and elementary flypes (or, maybe weighted elementary
flypes).
Third, understanding the connection between an interlocking homogeneous twisting
standard cabling torus and the Etnyre-Honda uniform thickness property
is a natural direction for future research.  The reader should realize that
construction of $\cR_{K(P,Q)}$ in \S\ref{section: building intlocking presentations}
can be varied by altering the sequencing choice of the sub-fiber discs and the
cabling/twisting illustrated in Figure \ref{figure:discs-plus-blocks-sequence}.
This spectrum of choices is conceivable a source of further examples of
transversal/Legendrian knots having the same classical invariants
but multiple transversal/Legendrian classes.

\end{document}